# Sufficient Conditions for the Exact Relaxation of Complementarity Constraints for Storages in Multi-period OPF Problems

Qi Wang, Wenchuan Wu, *Fellow, IEEE*, Chenhui Lin, Shuwei Xu, Xueliang Li

*Abstract*—Storage-concerned Optimal Power Flow (OPF) with complementarity constraints is highly non-convex and intractable. In this paper, we propose two generalized sufficient conditions which guarantee no simultaneous charging and discharging (SCD) in the relaxed multi-period OPF excluding the complementarity constraints. Moreover, we prove that the regions on the locational marginal prices (LMPs) formed by the proposed two conditions both contain the other existing representative ones. We also generalize the application premise of sufficient conditions from the positive electricity price requirements to the negative electricity price scenarios. In contrast to participating solely in the energy dispatch, we highlight that offering reserve services can broaden the application range of relaxation conditions for storages. The case studies verify the exactness and advantages of the proposed conditions.

## I. Introduction

### A. Background

Recently, there has been a rising trend in storage-plus-generation co-located hybrid resources, such as renewable energy paired with storages, and thermal units equipped with storages [1]. Furthermore, the Federal Energy Regulatory Commission (FERC) issued Order 841 in 2018 [2], which requires that each Regional Transmission Organization (RTO) and Independent System Operator (ISO) should develop a wholesale market participation model to facilitate the participation of storages in energy and ancillary services, through including storages' parameters that represent their distinct physical and operational characteristics.

In the conventional sense, the complementarity constraints of energy storage systems (ESSs) are introduced to avoid SCD, which render the whole optimization problem non-convex and intractable. To cope with this challenge, binary variables are usually employed to transform the complementarity constraints into mixed-integer programming (MIP). For instance, NYISO introduced its storage modeling within the dispatch model, incorporating binary variables for battery storages [3]. Besides, there is already a considerable body of literature exploring the involvement of energy storage in the joint energy-reserve optimization, most of which introduced integer variables to avoid SCD [4]-[7]. However, the MIP remains computationally complex as the number of ESSs or the whole time periods increases in multi-period OPF problems [8].

### B. Previous research

Recently, some literature has discussed the relaxation techniques for the complementarity constraints of ESSs, which can be roughly classified into four categories:

*i)* In this category, the complementarity constraints are reformulated as linear constraints about the charging/discharging active power and their upper bounds [9]-[11]. Nevertheless, the formulations still allow a considerable degree of simultaneous charging and discharging, that is, they cannot theoretically guarantee that SCD will be avoided.

*ii)* Feasible solution recovery strategies were introduced in [12], [13], with which the complementarity constraints could be met. Besides, the authors in [14] first supposed that charging and discharging efficiencies were both 1 to acquire a solution, which was probably not feasible in practice, and then readjusted the power of electric vehicles (EVs) for a new feasible solution. However, these methods can only ensure a suboptimal feasible solution with double computational burden.

*iii)* Based upon certain assumptions on the charging and discharging fees, some other relaxation methods were developed in [15]-[18]. But these methods may not be general enough to be applicable to all possible fees. Furthermore, a method to avoid the SCD by modifying the objective function was proposed in [19], but no prerequisite was given for its establishment.

*iv)* Some literature employ KKT conditions to analyze under which conditions the optimal solution satisfies the complementarity constraints [8], [20]-[24]. Wherein sufficient relaxation conditions were proposed for the economic dispatch (ED) problem of transmission networks in [20] and [21]. And an analysis was provided in [22] on the relaxation of grid-connected ESSs in a DCOPF problem for model predictive control. While the conditions presented in [22] not only rely on some strong assumptions on the objective function, but also are posterior, which limits their applications. As for ACOPF, the interaction between LMPs and the optimal operation of ESSs was studied in [23]. Besides, [8] and [24] provided preliminary analyses for the distribution network (DN) loss-minimizing problems, and derived sufficient conditions.

While the *condition C3* in [8] is rather restrictive, especially under the high penetration of renewable energy [8]. It's also noteworthy that the majority of existing relaxation conditions all require that the LMPs of the storage-connected buses should be

Manuscript received XX, 2023. This work was supported by the National Key Research and Development Plan of China under Grant 2022YFB2402900 and the S&T Program of State Grid Corporation of China under Grant 52060023001T (Key Techniques of Adaptive Grid Integration and Active Synchronization for Extremely High Penetration Distributed Photovoltaic Power Generation).

Q. Wang, W. Wu, C. Lin and S. Xu are with the Department of Electrical Engineering, Tsinghua University, Beijing 100084, China. (Corresponding Author: Wenchuan Wu, email: wuwench@tsinghua.edu.cn). X. Li is with State Grid Shandong Electric Power Co., Ltd., Jinan 250013, China.

non-negative, such as, the *relaxation condition of group C* in [21], and the *relaxation conditions* in [22]-[24]. Even though these conditions can be met in most of the cases, a growing number of negative electricity price scenarios appear with increasing renewable energy integration and transmission congestion. For example, [25] shows that the share of negative electricity prices is nearly 10% in February, 2022 in the Real-Time Market of ERCOT, as shown in Table I. In such cases, the aforementioned relaxation conditions may no longer apply. Considering that renewable energy and congestion are the major incentives for ESSs' deployment, the negative LMP scenarios can be expected to occur occasionally. Therefore, this issue need to be addressed.

In summary, a more general exact relaxation method should be investigated for multi-period OPF problems, which can be applied to a wider range of scenarios, such as negative LMPs.

TABLE I
PERCENTAGES OF NEGATIVE LMPs IN THE RTM OF ERCOT IN 2022 BY MONTH

| Month | LMP < 0 $/MWh | LMP < -10 $/MWh | LMP < -20 $/MWh | LMP < -30 $/MWh |
|---|---|---|---|---|
| 1 | 2.99% | 0.34% | 0.06% | 0.002% |
| 2 | 9.62% | 3.64% | 1.79% | 0.22% |
| 3 | 6.74% | 1.66% | 0.54% | 0.01% |
| 4 | 4.84% | 1.67% | 0.69% | 0.006% |
| 5 | 4.25% | 0.34% | 0.06% | 0.009% |
| 6 | 1.52% | 0.005% | 0.00% | 0.00% |
| 7 | 0.004% | 0.00% | 0.00% | 0.00% |
| 8 | 0.006% | 0.00% | 0.00% | 0.00% |
| 9 | 0.42% | 0.29% | 0.17% | 0.008% |
| 10 | 2.74% | 0.58% | 0.18% | 0.00% |
| 11 | 4.30% | 0.80% | 0.14% | 0.038% |
| 12 | 6.42% | 0.75% | 0.20% | 0.00% |

*C. Contributions*

This paper first presents a multi-period active and reactive power coordinated optimization model, considering energy and reserve services simultaneously. Subsequently, we conduct an in-depth analysis of the precise relaxation conditions for energy and reserve provision by energy storages. It's essential to emphasize that the proposed relaxation conditions can be applicable to both transmission and distribution grids. The contributions of this work are summarized as follows:

1) Two sufficient relaxation conditions are proposed that provably ensure no SCD in the relaxed multi-period OPF problems, such as storages participated ED, joint energy-reserve optimization and VAR optimization. And the application premise of relaxation conditions is generalized from positive electricity prices to negative LMP scenarios. Besides, the proposed two sufficient conditions do not have any requirement on the specific objective forms and power flow models, thus they hold for general cases. The accuracy and advantages of the proposed conditions are demonstrated in numerical tests.
2) To the best of our knowledge, we are the first to establish the relative inclusion relationships among multiple energy storage relaxation conditions. We prove that the regions about LMPs formulated by the proposed two conditions both contain the other existing ones in [20]-[24].
3) Prior relevant literature has solely focused on the relaxation conditions for storage participating in the energy dispatch [8]-[24], while disregarding its potential involvement in the reserve service. In light of this gap, we investigate the impact of storages involvement in reserve services on the relaxation conditions. Theoretical analysis indicates that offering reserve services can broaden the application range of relaxation conditions for storages, which is verified by case studies. As far as we know, such an effort has not been implemented in the previous studies.

The rest of this work is structured as follows: Section II introduces a multi-period coordinated optimization model for active and reactive power, taking into account both energy and reserve services. In section III, we first give two derived sufficient conditions, followed by the inclusion relationships among multiple relaxation conditions. Subsequently, we also elaborate on the applicability of the proposed relaxation conditions across various scenarios, such as the multi-period ED scenario. Case studies and comparative analyses are dedicated in section IV. Finally, section V summarizes the conclusions.

II. MATHEMATICAL MODEL FORMULATION

Consider a power system represented by a graph $\Pi(\mathcal{N}, \mathbb{E})$ with bus set $\mathcal{N} = \{1,...,n\}$, and line set $\mathbb{E} = \{i,j\} \in \mathcal{N} \times \mathcal{N}$. The time periods are denoted as $t = \{1,2,...,T\}$, with a scheduling interval $\Delta t$. It is crucial to emphasize that the given model formulation is not the innovation point of this work. Offering this comprehensive model is intended to facilitate the subsequent descriptions and proofs of the relaxation conditions.

*A. Objective function*

The objective function aims at seeking the optimal portfolio of active and reactive power to minimize the total production cost, which takes into account various components, including the generation cost, the reserve costs of conventional units and ESSs, penalties associated with renewable energy curtailment, energy storage arbitrage, and penalties for the simultaneous charging and discharging of storages.

$$\min \text{obj} = \sum_{t=1}^{T} \left\{ \sum_{n \in \Phi_{\text{ESS}}} \left[ g_{n,t}\left(p_{\text{ESS},n,t}^{\text{dc}}\right) + f_{n,t}\left(p_{\text{ESS},n,t}^{\text{ch}}\right) \right] \right.$$
$$+ \sum_{g \in \Phi_G} C_{g,t}\left(p_{g,t}^{\text{G}}\right) + \sum_{n \in \Phi_{\text{RG}}} C_{n,t}\left(p_{n,t}^{\text{RG}}\right)$$
$$+ \sum_{n \in \Phi_{\text{ESS}}} C_{n,t}^{\text{R}}\left(R_{n,t}^{\text{ch,ESS-}}, R_{n,t}^{\text{dc,ESS-}}, R_{n,t}^{\text{ch,ESS+}}, R_{n,t}^{\text{dc,ESS+}}\right)$$
$$\left. + \sum_{g \in \Phi_G} C_{g,t}^{\text{R}}\left(R_{g,t}^{\text{G+}}, R_{g,t}^{\text{G-}}\right) \right\} + C_{\text{ESS}}^{\text{Pen}} \quad (1)$$

$$C_{g,t}\left(p_{g,t}^{\text{G}}\right) = a_{2,g}\left(p_{g,t}^{\text{G}}\right)^2 + a_{1,g} p_{g,t}^{\text{G}} + a_{0,g} \quad (2)$$

$$C_{n,t}\left(p_{n,t}^{\text{RG}}\right) = b_{n,t} p_{n,t}^{\text{RG}} + \frac{\sigma_{\text{RG}}\left(p_{n,t}^{\text{RG}} - \widetilde{P}_{n,t}^{\text{RG}}\right)^2}{\widetilde{P}_{n,t}^{\text{RG}}} \quad (3)$$

$$C_{g,t}^{\text{R}}\left(R_{g,t}^{\text{G+}}, R_{g,t}^{\text{G-}}\right) = c_{g,t}^{\text{G+}} \cdot R_{g,t}^{\text{G+}} + c_{g,t}^{\text{G-}} \cdot R_{g,t}^{\text{G-}} \quad (4)$$

$$C_{n,t}^{\text{R}}\left(R_{n,t}^{\text{ch,ESS-}}, R_{n,t}^{\text{dc,ESS-}}, R_{n,t}^{\text{ch,ESS+}}, R_{n,t}^{\text{dc,ESS+}}\right)$$
$$= c_{n,t}^{\text{ESS+}} \cdot \left(R_{n,t}^{\text{ch,ESS+}} + R_{n,t}^{\text{dc,ESS+}}\right) + c_{n,t}^{\text{ESS-}} \cdot \left(R_{n,t}^{\text{ch,ESS-}} + R_{n,t}^{\text{dc,ESS-}}\right) \quad (5)$$

$$C_{\text{ESS}}^{\text{Pen}} = \sigma_{\text{ESS}} \sum_{t=1}^{T} \sum_{n \in \Phi_{\text{ESS}}} \left\{ p_{\text{ESS},n,t}^{\text{dc}} \left( \frac{1}{\eta_{\text{ESS},n}^{\text{dc}}} - 1 \right) + p_{\text{ESS},n,t}^{\text{ch}} \left( 1 - \eta_{\text{ESS},n}^{\text{ch}} \right) \right\} \quad (6)$$

where $g_{n,t}$ and $f_{n,t}$ are the discharging cost and charging fee of the $n^{\text{th}}$ ESS at time $t$. $C_{g,t}$ and $C_{n,t}$ denote the operational costs of the $g^{\text{th}}$ generator and the $n^{\text{th}}$ renewable generation (RG) at time $t$. $C_{g,t}^{R}$ and $C_{n,t}^{R}$ are the reserve costs of conventional units and storages. $C_{\text{ESS}}^{\text{Pen}}$ is the penalty term of storages to avoid the simultaneous charging and discharging paradox [19]. $p_{\text{ESS},n,t}^{\text{dc}}$ and $p_{\text{ESS},n,t}^{\text{ch}}$ are the discharging and charging active power of the $n^{\text{th}}$ ESS at time $t$. $p_{g,t}^{G}$ and $p_{n,t}^{\text{RG}}$ are the active power outputs of the $g^{\text{th}}$ generator and the $n^{\text{th}}$ RG at time $t$. $\bar{P}_{n,t}^{\text{RG}}$ is the predicted active power of the $n^{\text{th}}$ RG at time $t$. $R_{g,t}^{G+}$ and $R_{g,t}^{G-}$ denote the upward/downward reserve contributions of the $g^{\text{th}}$ generator at time $t$. $R_{n,t}^{\text{ch,ESS}-}, R_{n,t}^{\text{dc,ESS}-}, R_{n,t}^{\text{ch,ESS}+}, R_{n,t}^{\text{dc,ESS}+}$ represent the downward/upward reserves of the $n^{\text{th}}$ storage when it is in the charging/discharging status at time $t$, respectively. $a_{2,g}$, $a_{1,g}$ and $a_{0,g}$ are the coefficients of the quadratic, linear and constant terms in the operational cost of the $g^{\text{th}}$ generator. $b_{n,t}$ is the cost coefficient of the $n^{\text{th}}$ RG at time $t$. $c_{g,t}^{G+}$ and $c_{g,t}^{G-}$ denote the upward/downward reserve price of the $g^{\text{th}}$ generator at time $t$. $c_{n,t}^{\text{ESS}+}$ and $c_{n,t}^{\text{ESS}-}$ represent the upward/downward reserve price of the $n^{\text{th}}$ storage at time $t$. $\eta_{\text{ESS},n}^{\text{ch}}$ and $\eta_{\text{ESS},n}^{\text{dc}}$ are the charging and discharging efficiency of the $n^{\text{th}}$ ESS. $T$ is the whole time period. $\sigma_{\text{RG}}$ and $\sigma_{\text{ESS}}$ are the penalty coefficients. $\Phi_{\text{G}}, \Phi_{\text{RG}}$ and $\Phi_{\text{ESS}}$ are the index sets of generators, RGs and ESSs.

### B. Operational constraints

The operational constraints include not only power flow model, operational constraints of ESSs, generators, renewable generations and Static VAR compensators (SVCs), but also limits on voltage magnitudes, line transmission capacities and the minimum reserve requirement. Here, we denote the established model excluding the ESS' complementary constraint as the "*relaxed model*" of ACOPF.

#### 1) Power flow model

For the sake of generalizability, the AC power flow equations in polar coordinates are employed, ensuring the applicability of the proposed method in both transmission and distribution grids.

$$P_{ij,t} = \frac{1}{\tau_{ij,t}^2} g_{ij}^{\varepsilon} V_{i,t}^2 - \frac{1}{\tau_{ij,t}} V_{i,t} V_{j,t} \left[ g_{ij}^{\varepsilon} \cos(\theta_{i,t} - \theta_{j,t} - \phi_{ij,t}) + b_{ij}^{\varepsilon} \sin(\theta_{i,t} - \theta_{j,t} - \phi_{ij,t}) \right], \forall ij \in \mathbb{E}, \forall t \quad (7)$$

$$P_{ji,t} = g_{ij}^{\varepsilon} V_{j,t}^2 - \frac{1}{\tau_{ij,t}} V_{i,t} V_{j,t} \left[ g_{ij}^{\varepsilon} \cos(\theta_{j,t} - \theta_{i,t} + \phi_{ij,t}) + b_{ij}^{\varepsilon} \sin(\theta_{j,t} - \theta_{i,t} + \phi_{ij,t}) \right], \forall ij \in \mathbb{E}, \forall t \quad (8)$$

$$Q_{ij,t} = -\frac{1}{\tau_{ij,t}^2} \left( b_{ij}^{\varepsilon} + \frac{b_{ij}^{C}}{2} \right) V_{i,t}^2 - \frac{1}{\tau_{ij,t}} V_{i,t} V_{j,t} \left[ g_{ij}^{\varepsilon} \sin(\theta_{i,t} - \theta_{j,t} - \phi_{ij,t}) - b_{ij}^{\varepsilon} \cos(\theta_{i,t} - \theta_{j,t} - \phi_{ij,t}) \right], \forall ij \in \mathbb{E}, \forall t \quad (9)$$

$$Q_{ji,t} = -\left( b_{ij}^{\varepsilon} + \frac{b_{ij}^{C}}{2} \right) V_{j,t}^2 - \frac{1}{\tau_{ij,t}} V_{i,t} V_{j,t} \left[ g_{ij}^{\varepsilon} \sin(\theta_{j,t} - \theta_{i,t} + \phi_{ij,t}) - b_{ij}^{\varepsilon} \cos(\theta_{j,t} - \theta_{i,t} + \phi_{ij,t}) \right], \forall ij \in \mathbb{E}, \forall t \quad (10)$$

$$p_{j,t} - \sum_{i:ji \in \mathbb{E}} P_{ji,t} - \sum_{i:ij \in \mathbb{E}} P_{ji,t} - V_{j,t}^2 g_j^s = 0, \forall j \in \mathcal{N}, \forall t \quad \left( \lambda_{j,t}^{P} \right) \quad (11)$$

$$q_{j,t} - \sum_{i:ji \in \mathbb{E}} Q_{ji,t} - \sum_{i:ij \in \mathbb{E}} Q_{ji,t} + V_{j,t}^2 b_j^s = 0, \forall j \in \mathcal{N}, \forall t \quad (12)$$

where $P_{ij,t}$ and $Q_{ij,t}$ denote the active and reactive power flow of line $ij$ at time $t$. $V_{i,t}$ is the voltage magnitude of bus $i$ at time $t$. $\theta_{i,t}$ represents the phase angle of bus $i$ at time $t$. $\phi_{ij,t}$ is the transformer phase shift angle of branch $ij$ at time $t$. $p_{j,t}$ and $q_{j,t}$ are the net injected active and reactive power at bus $j$ and time $t$. $\tau_{ij,t}$ represents the transformer tap ratio of branch $ij$ at time $t$. $g_{ij}^{\varepsilon}$ and $b_{ij}^{\varepsilon}$ are the conductance and susceptance of branch $ij$. $b_{ij}^{C}$ denotes the line charging susceptance of branch $ij$. $g_j^s$ and $b_j^s$ are the shunt conductance and susceptance of bus $j$. $\lambda_{j,t}^{P}$ is the multiplier of the active power balance constraint at bus $j$ and time $t$.

$$p_{j,t} = \sum_{g \in \Phi_{\text{G},j}} p_{g,t}^{G} + \sum_{n \in \Phi_{\text{RG},j}} p_{n,t}^{\text{RG}} + \sum_{n \in \Phi_{\text{ESS},j}} \left( p_{\text{ESS},n,t}^{\text{dc}} - p_{\text{ESS},n,t}^{\text{ch}} \right) - p_{j,t}^{D} \quad (13)$$

$$q_{j,t} = \sum_{g \in \Phi_{\text{G},j}} q_{g,t}^{G} + \sum_{n \in \Phi_{\text{RG},j}} q_{n,t}^{\text{RG}} + \sum_{n \in \Phi_{\text{ESS},j}} q_{\text{ESS},n,t} + \sum_{n \in \Phi_{\text{SVC},j}} q_{n,t}^{\text{SVC}} - q_{j,t}^{D} \quad (14)$$

where $p_{g,t}^{G}$ and $q_{g,t}^{G}$ are the active and reactive power output of the $g^{\text{th}}$ generator at time $t$. $p_{n,t}^{\text{RG}}$ and $q_{n,t}^{\text{RG}}$ are the active and reactive power outputs of the $n^{\text{th}}$ RG at time $t$. $p_{\text{ESS},n,t}^{\text{dc}}$ and $p_{\text{ESS},n,t}^{\text{ch}}$ are the discharging and charging active power of the $n^{\text{th}}$ ESS at time $t$. $q_{\text{ESS},n,t}$ and $q_{n,t}^{\text{SVC}}$ denote the reactive power outputs of the $n^{\text{th}}$ ESS and the $n^{\text{th}}$ SVC at time $t$, respectively. $p_{j,t}^{D}$ and $q_{j,t}^{D}$ are the active and reactive power demands of the load at bus $j$ and time $t$. $\Phi_{\text{G},j}, \Phi_{\text{RG},j}, \Phi_{\text{ESS},j}$ and $\Phi_{\text{SVC},j}$ are the index set of generators, RGs, ESSs and SVCs connected to bus $j$.

#### 2) Operational constraints of ESSs

To take full advantage of the four-quadrant operation capabilities of ESSs, their operation should comply with the following constraints. The following constraints (17)-(21) together ensure the attainability of spinning reserves from storages by considering both the charging/discharging power and energy limits. The ESSs' SOC at the initial and final scheduling intervals are restricted to be consistent to guarantee a reliable energy planning for the next T-period dispatch cycle, as shown in (22).

$$-p_{\text{ESS},n,t}^{\text{ch}} \leq 0, \quad p_{\text{ESS},n,t}^{\text{ch}} \leq \bar{P}_{\text{ESS},n}^{\text{ch}}, \quad \left( \lambda_{\text{ESS},n,t}^{\text{ch},1}, \lambda_{\text{ESS},n,t}^{\text{ch},2} \right) \quad (15)$$

$$-p_{\text{ESS},n,t}^{\text{dc}} \leq 0, \quad p_{\text{ESS},n,t}^{\text{dc}} \leq \bar{P}_{\text{ESS},n}^{\text{dc}}, \quad \left( \lambda_{\text{ESS},n,t}^{\text{dc},1}, \lambda_{\text{ESS},n,t}^{\text{dc},2} \right) \quad (16)$$

ESSs can provide upward reserves by either reducing the charging power or increasing the discharging power. Similarly, downward reserve capacity can be offered by charging more or discharging less from their predetermined power schedules [26].

$$p_{\text{ESS},n,t}^{\text{ch}} + R_{n,t}^{\text{ch,ESS}-} \leq \overline{P}_{\text{ESS},n}^{\text{ch}}, \quad \left(\lambda_{\text{ESS},n,t}^{\text{RD},1}\right) \tag{17}$$

$$R_{n,t}^{\text{dc,ESS}-} \leq p_{\text{ESS},n,t}^{\text{dc}}, \quad \left(\lambda_{\text{ESS},n,t}^{\text{RD},2}\right)$$

$$p_{\text{ESS},n,t}^{\text{dc}} + R_{n,t}^{\text{dc,ESS}+} \leq \overline{P}_{\text{ESS},n}^{\text{dc}}, \quad \left(\lambda_{\text{ESS},n,t}^{\text{RU},1}\right) \tag{18}$$

$$R_{n,t}^{\text{ch,ESS}+} \leq p_{\text{ESS},n,t}^{\text{ch}}, \quad \left(\lambda_{\text{ESS},n,t}^{\text{RU},2}\right)$$

$$R_{n,t}^{\text{ch,ESS}-}, R_{n,t}^{\text{dc,ESS}-}, R_{n,t}^{\text{ch,ESS}+}, R_{n,t}^{\text{dc,ESS}+} \geq 0 \tag{19}$$

$$\left(1-\delta_{\text{ESS},n}\right)^t E_{\text{ESS},n,0} + \sum_{\tau=1}^{t}\left(1-\delta_{\text{ESS},n}\right)^{t-\tau}$$
$$\times\left(\left(p_{\text{ESS},n,\tau}^{\text{ch}} - R_{n,\tau}^{\text{ch,ESS}+}\right)\cdot\eta_{\text{ESS},n}^{\text{ch}} - \left(p_{\text{ESS},n,\tau}^{\text{dc}} + R_{n,\tau}^{\text{dc,ESS}+}\right)/\eta_{\text{ESS},n}^{\text{dc}}\right)\Delta t \geq \underline{E}_{\text{ESS},n}, \quad \left(\lambda_{\text{ESS},n,t}^{\text{SOC},1}\right) \tag{20}$$

$$\left(1-\delta_{\text{ESS},n}\right)^t E_{\text{ESS},n,0} + \sum_{\tau=1}^{t}\left(1-\delta_{\text{ESS},n}\right)^{t-\tau}$$
$$\times\left(\left(p_{\text{ESS},n,\tau}^{\text{ch}} + R_{n,\tau}^{\text{ch,ESS}-}\right)\cdot\eta_{\text{ESS},n}^{\text{ch}} - \left(p_{\text{ESS},n,\tau}^{\text{dc}} - R_{n,\tau}^{\text{dc,ESS}-}\right)/\eta_{\text{ESS},n}^{\text{dc}}\right)\Delta t \leq \overline{E}_{\text{ESS},n}, \quad \left(\lambda_{\text{ESS},n,t}^{\text{SOC},2}\right) \tag{21}$$

$$E_{\text{ESS},n,T} = \left(1-\delta_{\text{ESS},n}\right)^T E_{\text{ESS},n,0} + \sum_{\tau=1}^{T}\left(1-\delta_{\text{ESS},n}\right)^{T-\tau}$$
$$\times\left(p_{\text{ESS},n,\tau}^{\text{ch}}\cdot\eta_{\text{ESS},n}^{\text{ch}} - p_{\text{ESS},n,\tau}^{\text{dc}}/\eta_{\text{ESS},n}^{\text{dc}}\right)\Delta t = E_{\text{ESS},n,0} \tag{22}$$

$$\left(p_{\text{ESS},n,t}^{\text{dc}}\right)^2 + \left(q_{\text{ESS},n,t}\right)^2 \leq \left(S_{\text{ESS},n}\right)^2, \quad \left(\lambda_{\text{ESS},n,t}^{\text{S},1}\right) \tag{23}$$

$$\left(p_{\text{ESS},n,t}^{\text{ch}}\right)^2 + \left(q_{\text{ESS},n,t}\right)^2 \leq \left(S_{\text{ESS},n}\right)^2, \quad \left(\lambda_{\text{ESS},n,t}^{\text{S},2}\right) \tag{24}$$

where $E_{\text{ESS},n,t}$ represents the SOC level of the $n^{\text{th}}$ ESS at time $t$. $E_{\text{ESS},n,0}$ denotes the initial SOC of the $n^{\text{th}}$ ESS. $\delta_{\text{ESS},n}$ is the self-discharge rate of the $n^{\text{th}}$ ESS. $\eta_{\text{ESS},n}^{\text{ch}}$ and $\eta_{\text{ESS},n}^{\text{dc}}$ are the charging and discharging efficiency of the $n^{\text{th}}$ ESS. $\underline{E}_{\text{ESS},n}$ are $\overline{E}_{\text{ESS},n}$ the lower and upper bounds of the SOC level of the $n^{\text{th}}$ ESS, respectively. $\overline{P}_{\text{ESS},n}^{\text{ch}}$ and $\overline{P}_{\text{ESS},n}^{\text{dc}}$ are the upper limits of charging and discharging active power of the $n^{\text{th}}$ ESS, respectively. $S_{\text{ESS},n}$ is the apparent power of the $n^{\text{th}}$ ESS. $\lambda_{\text{ESS},n,t}^{\text{ch},1}$, $\lambda_{\text{ESS},n,t}^{\text{ch},2}$, $\lambda_{\text{ESS},n,t}^{\text{dc},1}$, $\lambda_{\text{ESS},n,t}^{\text{dc},2}$, $\lambda_{\text{ESS},n,t}^{\text{RD},1}$, $\lambda_{\text{ESS},n,t}^{\text{RD},2}$, $\lambda_{\text{ESS},n,t}^{\text{RU},1}$, $\lambda_{\text{ESS},n,t}^{\text{RU},2}$, $\lambda_{\text{ESS},n,t}^{\text{SOC},1}$, $\lambda_{\text{ESS},n,t}^{\text{SOC},2}$, $\lambda_{\text{ESS},n,t}^{\text{S},1}$ and $\lambda_{\text{ESS},n,t}^{\text{S},2}$ are the multipliers of the corresponding constraints.

Besides, we can employ (25) as an enhanced complement to the aforementioned ESS model.

$$\frac{p_{\text{ESS},n,t}^{\text{ch}}}{\overline{P}_{\text{ESS},n}^{\text{ch}}} + \frac{p_{\text{ESS},n,t}^{\text{dc}}}{\overline{P}_{\text{ESS},n}^{\text{dc}}} \leq 1, \forall n \in \Phi_{\text{ESS}}, \forall t, \quad \left(\lambda_{\text{ESS},n,t}^{\text{relax}}\right) \tag{25}$$

*3) Other operational constraints*

The operational constraints of generators are presented as:

$$\underline{P}_g^{\text{G}} \leq p_{g,t}^{\text{G}} \leq \overline{P}_g^{\text{G}}, \forall g \in \Phi_{\text{G}} \tag{26}$$

$$\underline{Q}_g^{\text{G}} \leq q_{g,t}^{\text{G}} \leq \overline{Q}_g^{\text{G}}, \forall g \in \Phi_{\text{G}} \tag{27}$$

$$-RD_g \Delta t \leq p_{g,t+1}^{\text{G}} - p_{g,t}^{\text{G}} \leq RU_g \Delta t \tag{28}$$

$$0 \leq R_{g,t}^{\text{G}+} \leq RU_g \Delta t,$$
$$R_{g,t}^{\text{G}+} \leq \overline{P}_g^{\text{G}} - p_{g,t}^{\text{G}}, \forall g \in \Phi_{\text{G}} \tag{29}$$

$$0 \leq R_{g,t}^{\text{G}-} \leq RD_g \Delta t,$$
$$R_{g,t}^{\text{G}-} \leq p_{g,t}^{\text{G}} - \underline{P}_g^{\text{G}}, \forall g \in \Phi_{\text{G}} \tag{30}$$

$$\sum_{g \in \Phi_{\text{G}}} R_{g,t}^{\text{G}+} + \sum_{n \in \Phi_{\text{ESS}}}\left(R_{n,t}^{\text{ch,ESS}+} + R_{n,t}^{\text{dc,ESS}+}\right) \geq \text{SRU}_t,$$
$$\sum_{g \in \Phi_{\text{G}}} R_{g,t}^{\text{G}-} + \sum_{n \in \Phi_{\text{ESS}}}\left(R_{n,t}^{\text{ch,ESS}-} + R_{n,t}^{\text{dc,ESS}-}\right) \geq \text{SRD}_t \tag{31}$$

where $\underline{P}_g^{\text{G}}$ and $\overline{P}_g^{\text{G}}$ are the lower and upper bounds of the active power output of the $g^{\text{th}}$ generator. $\underline{Q}_g^{\text{G}}$ and $\overline{Q}_g^{\text{G}}$ represent the lower and upper bounds of the reactive output of the $g^{\text{th}}$ generator. $RU_g$ and $RD_g$ denote the upward and downward ramp rates of the $g^{\text{th}}$ generator. $\text{SRU}_t$ and $\text{SRD}_t$ are the system-wide upward and downward spinning reserve capacity requirements, which can be determined based on the historical uncertainty data of renewable energy and loads [5].

The operational constraints of RGs are given as follows:

$$\underline{P}_{n,t}^{\text{RG}} \leq p_{n,t}^{\text{RG}} \leq \tilde{P}_{n,t}^{\text{RG}}, \forall n \in \Phi_{\text{RG}} \tag{32}$$

$$\left(p_{n,t}^{\text{RG}}\right)^2 + \left(q_{n,t}^{\text{RG}}\right)^2 \leq \left(S_n^{\text{RG}}\right)^2, \forall n \in \Phi_{\text{RG}} \tag{33}$$

where $\tilde{P}_{n,t}^{\text{RG}}$ is the predicted active power of the $n^{\text{th}}$ RG at time $t$. $\underline{P}_{n,t}^{\text{RG}}$ is the lower bound of the active power output of the $n^{\text{th}}$ RG at time $t$. And $S_n^{\text{RG}}$ is the apparent power of the $n^{\text{th}}$ RG.

The operational constraints of SVCs are described as:

$$\underline{Q}_n^{\text{SVC}} \leq q_{n,t}^{\text{SVC}} \leq \overline{Q}_n^{\text{SVC}}, \forall n \in \Phi_{\text{SVC}} \tag{34}$$

where $\underline{Q}_n^{\text{SVC}}$ and $\overline{Q}_n^{\text{SVC}}$ are the lower and upper bounds of the reactive power outputs of the $n^{\text{th}}$ SVC, respectively.

The system operation also needs to meet the voltage security constraints and line transmission capacity limits.

$$\underline{V}_i \leq V_{i,t} \leq \overline{V}_i, \forall i \in \mathcal{N} \tag{35}$$

$$P_{ij,t}^2 + Q_{ij,t}^2 \leq \overline{S}_{ij}^2, \forall ij \in \mathcal{E} \tag{36}$$

where $\underline{V}_i$ and $\overline{V}_i$ are the lower and upper bounds of voltage magnitude of bus $i$, respectively. $\overline{S}_{ij}$ is the transmission capacity of branch $ij$.

### III. EXACT RELAXATION CONDITIONS FOR COMPLEMENTARITY CONSTRAINTS OF STORAGES

*A. Two sufficient conditions for exactly relaxing ESSs' complementary constraints*

**Theorem 1**: *The proposed two sufficient conditions to avoid SCD adopting the relaxed model are given as follows. This means that as long as any one of the following conditions is satisfied, the relaxation is guaranteed to be exact.*

*C1*: $\lambda_{j,t}^{\text{p}} > \lambda_{j,t}^{\text{p,C1}}$ \quad (37)

*C2*: $\lambda_{j,t}^{\text{p}} > \lambda_{j,t}^{\text{p,C2}}$ \quad (38)

where

$$\lambda_{j,t}^{\text{p,C1}} = \left[ \frac{\partial \text{obj}}{\partial \left(p_{\text{ESS},n,t}^{\text{ch}}\right)} / \eta_{\text{ESS},n}^{\text{dc}} + \frac{\partial \text{obj}}{\partial \left(p_{\text{ESS},n,t}^{\text{dc}}\right)} \cdot \eta_{\text{ESS},n}^{\text{ch}} \right] / \left( \eta_{\text{ESS},n}^{\text{ch}} - \frac{1}{\eta_{\text{ESS},n}^{\text{dc}}} \right)$$

$$- \left\{ \left[ \lambda_{\text{ESS},n,t}^{\text{ch,2}} + 2\lambda_{\text{ESS},n,t}^{\text{S,2}} \cdot p_{\text{ESS},n,t}^{\text{ch}} + \frac{\lambda_{\text{ESS},n,t}^{\text{relax}}}{\overline{P}_{\text{ESS},n}^{\text{ch}}} + \lambda_{\text{ESS},n,t}^{\text{RD,1}} \right] / \eta_{\text{ESS},n}^{\text{dc}} \right.$$

$$\left. + \left[ \lambda_{\text{ESS},n,t}^{\text{dc,2}} + 2\lambda_{\text{ESS},n,t}^{\text{S,1}} \cdot p_{\text{ESS},n,t}^{\text{dc}} + \frac{\lambda_{\text{ESS},n,t}^{\text{relax}}}{\overline{P}_{\text{ESS},n}^{\text{dc}}} + \lambda_{\text{ESS},n,t}^{\text{RU,1}} \right] \cdot \eta_{\text{ESS},n}^{\text{ch}} \right\} / \left( \frac{1}{\eta_{\text{ESS},n}^{\text{dc}}} - \eta_{\text{ESS},n}^{\text{ch}} \right)$$
(39)

$$\lambda_{j,t}^{\text{p,C2}} = \left[ \frac{\partial \text{obj}}{\partial \left(p_{\text{ESS},n,t}^{\text{ch}}\right)} / \eta_{\text{ESS},n}^{\text{dc}} + \frac{\partial \text{obj}}{\partial \left(p_{\text{ESS},n,t}^{\text{dc}}\right)} \cdot \eta_{\text{ESS},n}^{\text{ch}} \right] / \left( \eta_{\text{ESS},n}^{\text{ch}} - \frac{1}{\eta_{\text{ESS},n}^{\text{dc}}} \right) \quad (40)$$

where $\lambda_{j,t}^{\text{p}}$ is the multiplier of the active power balance constraint at bus $j$ and time $t$. $\lambda_{\text{ESS},n,t}^{\text{ch,2}}$, $\lambda_{\text{ESS},n,t}^{\text{dc,2}}$, $\lambda_{\text{ESS},n,t}^{\text{S,1}}$, $\lambda_{\text{ESS},n,t}^{\text{S,2}}$, $\lambda_{\text{ESS},n,t}^{\text{RD,1}}$, $\lambda_{\text{ESS},n,t}^{\text{RU,1}}$ and $\lambda_{\text{ESS},n,t}^{\text{relax}}$ are multipliers of the corresponding constraints.

*Proof*: The detailed proof is given in *Appendix A* in the supplementary file [27].
∎

*Remark 1*: Considering that the **Condition C1** in *Theorem 1* is posterior, that is, dependent on the results. Thus, we further propose the **Condition C2**. When the terms in the objective function regarding the charging and discharging active power of ESSs are all linear, the right-hand term of formula (40) becomes a constant, which facilitates the prior verification. That is, **Condition C2** can be checked before the problem solving since the LMP can be forecasted based on the historical data, which facilitates its application. The practical implementation procedure of the derived sufficient condition C2 will be detailed in Section IV.

*Remark 2*: It should also be emphasized that the two proposed relaxation conditions are independent from the forms of objective functions and power flow models, and thus the conclusions are general. In other words, when employing other power flow models [28]-[29], the derived relaxation conditions remain applicable.

### B. The inclusion relationships among multiple relaxation conditions

Firstly, the relationship between *Condition C1* and *Condition C2* in *Theorem 1* is described in *Lemma 1*.

*Lemma 1*: The region on LMPs formulated by **Condition C1** contains that of **Condition C2**.

*Proof*: The detailed proof is given in *Appendix B* in the supplementary file [27].
∎

Subsequently, we present the following lemma to discuss the inclusion relationships between the proposed two relaxation conditions and the other existing ones in [20]-[24]. Although the conditions given in [20]-[22] are based upon the ED model, similar relaxation conditions can be derived for ACOPF, which are introduced in *Appendix C* in the supplementary file [27].

*Lemma 2*: The regions about LMPs formulated by the proposed sufficient conditions in **Theorem 1** both contain those of the other six existing relaxation conditions:
 i) the Relaxation conditions in [20] (i.e., the Relaxation conditions of group A in [21]);
 ii) the Relaxation conditions of group B in [21];
 iii) the Relaxation conditions of group C in [21];
 iv) the Relaxation conditions in [22];
 v) the Relaxation conditions in [23];
 vi) the Relaxation conditions in [24].

*Proof*: The detailed proofs are given in *Appendix D* in the supplementary file [27]. And elaborate descriptions of these relaxation conditions are given in *Appendixes C and D* to facilitate the proofs.
∎

*Remark 3*: According to the proofs in *Appendix D*, we can conclude that the presented **Lemma 2** does not rely on the specific objective forms and power flow models, and thus it is generally valid. Moreover, Fig. 1 illustrates the inclusion relationships between different relaxation conditions given in **Lemma 2**. As far as we know, it is the first time that the inclusion relationships among multiple energy storage relaxation conditions are established.

### C. Elaborating on the applicability of the proposed relaxation conditions across various scenarios

#### 1) Multi-period ED scenario

It's essential to clarify that the novelty of the proposed relaxation conditions does not arise from the transition of the scenarios from ED to ACOPF, but rather from our methodological innovations, aiming to achieve broader relaxation conditions. In other words, the generalization comes from a novel demonstration of the conditions rather than the ACOPF formulation.

For a clearer demonstration, we will analyze the relaxation conditions for ESSs *under multi-period ED scenarios*. Similar to the aforementioned *Conditions C1* and *C2*, the following two relaxation conditions are introduced. And please refer to *Appendix E* in the supplementary file [27] for detailed proofs.

**C1$^{\text{ED}}$**: $\text{LMP}_{j,t} > \lambda_{j,t}^{\text{ED,C1}}$ (41)

**C2 $^{\text{ED}}$**: $\text{LMP}_{j,t} > \lambda_{j,t}^{\text{ED,C2}}$ (42)

where $\text{LMP}_{j,t} = \lambda_{j,t}^{\text{p}} + \sum_{j} SF_{j-i} \left( \lambda_{\text{L},j}^{1}(t) - \lambda_{\text{L},j}^{2}(t) \right)$. $SF_{j-i}$ denotes the shift factor for bus $i$ on line $j$. $\lambda_{\text{L},j}^{2}(t)$ and $\lambda_{\text{L},j}^{1}(t)$ represent the multipliers corresponding to the upper and lower limits of line capacity characterized by shift factors.

$$\lambda_{j,t}^{\text{ED,C1}} = \left[ \frac{\partial \text{obj}}{\partial \left(p_{\text{ESS},n,t}^{\text{ch}}\right)} / \eta_{\text{ESS},n}^{\text{dc}} + \frac{\partial \text{obj}}{\partial \left(p_{\text{ESS},n,t}^{\text{dc}}\right)} \cdot \eta_{\text{ESS},n}^{\text{ch}} \right] / \left( \eta_{\text{ESS},n}^{\text{ch}} - \frac{1}{\eta_{\text{ESS},n}^{\text{dc}}} \right)$$

$$- \left\{ \left[ \lambda_{\text{ESS},n,t}^{\text{ch,2}} + \frac{\lambda_{\text{ESS},n,t}^{\text{relax}}}{\overline{P}_{\text{ESS},n}^{\text{ch}}} + \lambda_{\text{ESS},n,t}^{\text{RD,1}} \right] / \eta_{\text{ESS},n}^{\text{dc}} \right.$$ (43)

$$\left. + \left[ \lambda_{\text{ESS},n,t}^{\text{dc,2}} + \frac{\lambda_{\text{ESS},n,t}^{\text{relax}}}{\overline{P}_{\text{ESS},n}^{\text{dc}}} + \lambda_{\text{ESS},n,t}^{\text{RU,1}} \right] \cdot \eta_{\text{ESS},n}^{\text{ch}} \right\} / \left( \frac{1}{\eta_{\text{ESS},n}^{\text{dc}}} - \eta_{\text{ESS},n}^{\text{ch}} \right)$$

$$\lambda_{j,t}^{\text{ED,C2}} = \left[ \frac{\partial \text{obj}}{\partial \left(p_{\text{ESS},n,t}^{\text{ch}}\right)} / \eta_{\text{ESS},n}^{\text{dc}} + \frac{\partial \text{obj}}{\partial \left(p_{\text{ESS},n,t}^{\text{dc}}\right)} \cdot \eta_{\text{ESS},n}^{\text{ch}} \right] / \left( \eta_{\text{ESS},n}^{\text{ch}} - \frac{1}{\eta_{\text{ESS},n}^{\text{dc}}} \right) \quad (44)$$

Through comparisons, it can be observed that *Condition C2$^{\text{ED}}$* is consistent with *Condition C2*. And in comparison to *Condition C1*, *Condition C1$^{\text{ED}}$* only omits the terms related to

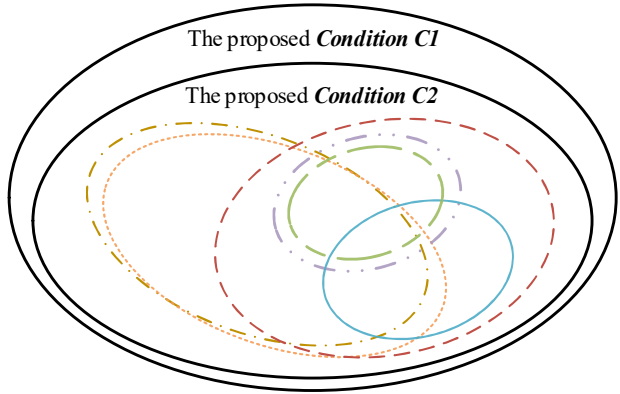

Fig. 1. The inclusion relationships between different relaxation conditions.

i) the Relaxation conditions in [20] (i.e. the Relaxation conditions of group A in [21])
ii) the Relaxation conditions of group B in [21]
iii) the Relaxation conditions of group C in [21]
iv) the Relaxation conditions in [22]
v) the Relaxation conditions in [23]
vi) the Relaxation conditions in [24]

reactive power.

*2) The ACOPF scenario in which storages do not participate in reserve services*

When energy storages are not involvement in reserve services, it is equivalent to setting $R_{n,t}^{\text{ch,ESS}-}, R_{n,t}^{\text{dc,ESS}-}, R_{n,t}^{\text{ch,ESS}+}, R_{n,t}^{\text{dc,ESS}+}$ to zero in formulas (17)-(21). In such scenario, constraints (17)-(19) degrade to become consistent with (15)-(16). Correspondingly, the two relaxation conditions provided in *Theorem 1* are transformed into:

$$C1^{\text{NoR}}: \lambda_{j,t}^{\text{p}} > \lambda_{j,t}^{\text{NoR,C1}} \tag{45}$$

$$C2^{\text{NoR}}: \lambda_{j,t}^{\text{p}} > \lambda_{j,t}^{\text{NoR,C2}} \tag{46}$$

where

$$\lambda_{j,t}^{\text{NoR,C1}} = \left[ \frac{\partial \text{obj}}{\partial(p_{\text{ESS},n,t}^{\text{ch}})} / \eta_{\text{ESS},n}^{\text{dc}} + \frac{\partial \text{obj}}{\partial(p_{\text{ESS},n,t}^{\text{dc}})} \cdot \eta_{\text{ESS},n}^{\text{ch}} \right] / \left( \eta_{\text{ESS},n}^{\text{ch}} - \frac{1}{\eta_{\text{ESS},n}^{\text{dc}}} \right)$$

$$- \left\{ \left[ \lambda_{\text{ESS},n,t}^{\text{ch,2}} + 2\lambda_{\text{ESS},n,t}^{\text{S,2}} \cdot p_{\text{ESS},n,t}^{\text{ch}} + \frac{\lambda_{\text{ESS},n,t}^{\text{relax}}}{\overline{P}_{\text{ESS},n}^{\text{ch}}} \right] / \eta_{\text{ESS},n}^{\text{dc}} \right. \tag{47}$$

$$\left. + \left[ \lambda_{\text{ESS},n,t}^{\text{dc,2}} + 2\lambda_{\text{ESS},n,t}^{\text{S,1}} \cdot p_{\text{ESS},n,t}^{\text{dc}} + \frac{\lambda_{\text{ESS},n,t}^{\text{relax}}}{\overline{P}_{\text{ESS},n}^{\text{dc}}} \right] \cdot \eta_{\text{ESS},n}^{\text{ch}} \right\} / \left( \frac{1}{\eta_{\text{ESS},n}^{\text{dc}}} - \eta_{\text{ESS},n}^{\text{ch}} \right)$$

$$\lambda_{j,t}^{\text{NoR,C2}} = \left[ \frac{\partial \text{obj}}{\partial(p_{\text{ESS},n,t}^{\text{ch}})} / \eta_{\text{ESS},n}^{\text{dc}} + \frac{\partial \text{obj}}{\partial(p_{\text{ESS},n,t}^{\text{dc}})} \cdot \eta_{\text{ESS},n}^{\text{ch}} \right] / \left( \eta_{\text{ESS},n}^{\text{ch}} - \frac{1}{\eta_{\text{ESS},n}^{\text{dc}}} \right) \tag{48}$$

After making comparisons, it is apparent that *Condition C2*[NoR] is exactly the same as *Condition C2*. And *Condition C1*[NoR] differs from *Condition C1* solely in that it removes the terms related to the energy storage reserve constraints (17)-(18).

In summary, *Condition C2* is generally valid in various scenarios based on the above analyses. As for *Condition C1*, it remains applicable in multiple scenarios with the necessary adjustments to adapt to different models. The generality of the proposed relaxation conditions presented in *Theorem 1*, enhances their potentials for broader applications.

***Remark 4***: Subsequently, we investigate the impact of storage involvement in reserve services on the relaxation conditions.

***Lemma 3***: Offering reserve services has the potential to expand the application range of relaxation conditions for storages.

***Proof***: Recall that *Condition C1*[NoR] differs from *Condition C1* solely in its omittance of the terms related to the energy storage reserve constraints (17)-(18). Considering that the dual multipliers corresponding to (17)-(18) are all non-negative, thus, we can obtain that $\lambda_{j,t}^{\text{NoR,C1}} > \lambda_{j,t}^{\text{p,C1}}$. Combining that *Condition C2*[NoR] is consistent with *Condition C2*, it follows that $\lambda_{j,t}^{\text{NoR,C2}} = \lambda_{j,t}^{\text{p,C2}}$. In summary, we can conclude that providing reserve services holds the potential to broaden the applicability of the relaxation conditions for storages.

## IV. NUMERICAL TESTS

A 24-hour ACOPF problem is tested using the IEEE 69 bus system, including thirty-five ESSs, twenty-nine RGs, five generators and seven SVCs. The charging/discharging efficiencies of storages are randomly generated within the range of [0.9, 0.93]. Fig. 2 illustrates the total loads and available renewable energy. The detailed parameters about loads, ESSs, RGs, generators and SVCs are given in [30]. Moreover, we set self-discharge rate as 0.002, $\sigma_{\text{ESS}} = 3$, $\Delta t = 1\,\text{h}$, $c_{g,t}^{\text{G}+} = c_{g,t}^{\text{G}-} = 10\,\$/\text{MWh}$, $c_{n,t}^{\text{ESS}+} = c_{n,t}^{\text{ESS}-} = 5\,\$/\text{MWh}$. The solver for continuous optimization problems is IPOPT [31], and BONMIN [32] for MIP. Case studies are conducted in Matlab 2019b on a laptop with an Intel i7-10875H CPU and 24GB RAM.

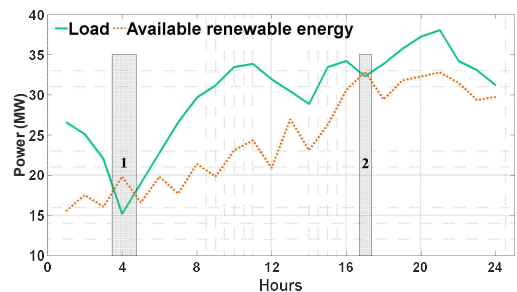

Fig. 2. Illustration of the total loads and available renewable energy.

### A. Practical application scheme

Researchers have been consistently investigating the LMP forecasting and have reported several effective approaches [33]-[35], such as deep learning neural network and generative adversarial network (GAN). To better illustrate the prediction accuracy, some test results for the day-ahead LMP prediction in [35] are displayed in Fig. 3. It is noticeable that the GAN-based approach can accurately capture distinct daily price patterns for different locations. As LMP forecasts become more accurate, it also facilitates the practical application of the proposed relaxation conditions.

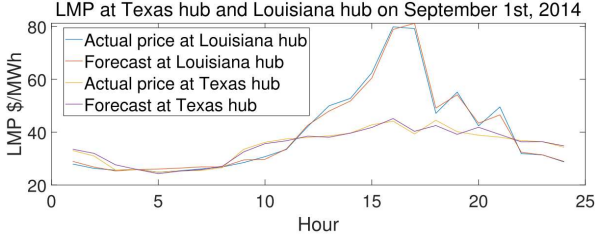

Fig. 3. The LMP prediction results reported in [35].

Based on the forecasted LMPs, it is possible to utilize the lower bound of the interval into which the actual LMP may fall, denoted as $\underline{\text{LMP}}_{j,t}$, to determine the exactness of the aforementioned sufficient conditions in real-world operational scenarios. More precisely, if
$$\left[ \frac{\partial \text{obj}}{\partial \left( p_{\text{ESS},n,t}^{\text{ch}} \right)} / \eta_{\text{ESS},n}^{\text{dc}} + \frac{\partial \text{obj}}{\partial \left( p_{\text{ESS},n,t}^{\text{dc}} \right)} \cdot \eta_{\text{ESS},n}^{\text{ch}} \right] / \left( \eta_{\text{ESS},n}^{\text{ch}} - \frac{1}{\eta_{\text{ESS},n}^{\text{dc}}} \right)$$
is significantly lower than $\underline{\text{LMP}}_{j,t}$, it can be deduced that *Condition C2* is highly likely to be met.

As a matter of fact, employing $\underline{\text{LMP}}_{j,t}$ rather than $\text{LMP}_{j,t}$, facilitates the applicability of the exact relaxation conditions. This is because forecasting the range within which the LMPs are likely to lie, is more easily implementable than precisely determining the actual LMP values in practical scenarios. Specifically, we can determine the lower bound $\underline{\text{LMP}}_{j,t}$ using the mean absolute percentage error (MAPE) as an approximate estimation of the standard deviation in LMP forecasting:

$$\underline{\text{LMP}}_{j,t} = \begin{cases} \left(1 - 5 \times \text{MAPE}\right) \times \text{LMP}_{j,t}^{\text{pred}}, & \text{if } \text{LMP}_{j,t}^{\text{pred}} > 0 \\ \left(1 + 5 \times \text{MAPE}\right) \times \text{LMP}_{j,t}^{\text{pred}}, & \text{otherwise} \end{cases} \quad (49)$$

where $\text{LMP}_{j,t}^{\text{pred}}$ is the predicted LMP.

In light of the insights derived from the above discussions, a practical scheme for using the exact sufficient *Condition C2* can be performed as follows. When *Condition C2* is not satisfied, an increase in the penalty coefficient $\sigma_{\text{ESS}}$, if appropriately adjusted, can guarantee the exactness of the relaxation. Theoretically, as long as the penalty coefficient $\sigma_{\text{ESS}}$ is large enough, *Condition C2* is guaranteed to be met, ensuring the exactness of the relaxation.

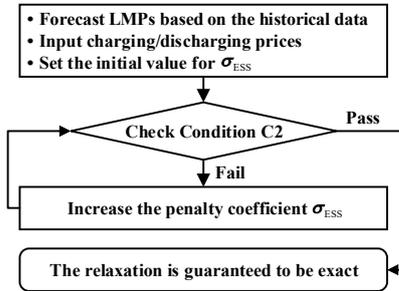

Fig. 4. Flowchart of the practical implementation of the derived sufficient condition C2.

### B. Day-ahead cases under negative LMP scenarios

In the following case studies, $f_{n,t}\left(p_{\text{ESS},n,t}^{\text{ch}}\right)$ and $g_{n,t}\left(p_{\text{ESS},n,t}^{\text{dc}}\right)$ are reformulated as: $f_{n,t}\left(p_{\text{ESS},n,t}^{\text{ch}}\right) = a_t^{\text{ch}} \cdot p_{\text{ESS},n,t}^{\text{ch}}$ and $g_{n,t}\left(p_{\text{ESS},n,t}^{\text{dc}}\right) = a_t^{\text{dc}} \cdot p_{\text{ESS},n,t}^{\text{dc}}$, where $a_t^{\text{ch}}$ and $a_t^{\text{dc}}$ are the charging and discharging prices at time $t$. Given that $a_t^{\text{dc}} = a_t^{\text{ch}} = 15$ \$/MWh and the renewable generation cost is -100 \$/MWh, representing scenarios with notably high penetrations of renewable energy, the optimization results are depicted in Fig. 5. In subfigure (a), the gray areas indicate that storage charging takes place during periods when the available renewable energy exceeds the load demand, as illustrated in Fig. 2. This is because storages absorbing surplus power helps minimize the renewable energy curtailment.

As shown in subfigure (b), the results not only demonstrate that both *Conditions C1* and *C2* hold in this scenario, but also verify that *Lemma 1* is satisfied. Furthermore, the negative LMPs appear exactly during the periods when the available renewable energy exceeds the load demand, that is, the gray areas. Correspondingly, it can be observed that there is no SCD in the optimization results of subfigure (a), verifying the exactness of the proposed conditions. In contrast, the relaxation conditions of Group C in [21] and conditions in [22]-[24] do not hold due to the existence of negative LMPs. As for the relaxation conditions of groups A and B in [21] and conditions in [20], none of them hold since the negative price values are quite small, illustrating the advantages of the proposed conditions under negative prices.

In subfigure (c), we compare relaxation conditions' boundaries with and without providing reserve services. It can be observed that offering reserve services has the potential to expand the application range of relaxation conditions for storages, which is consistent with *Lemma 3*, as presented by the yellow solid line and orange dotted line. Recall that *Condition C2$^{\text{NoR}}$* is consistent with *Condition C2*, thus, we do not make a distinction on *Condition C2* for these two scenarios.

In addition, subfigure (d) presents the upward/downward reserves supplied by the storage at bus 10. Take the 2nd grey area as an example, it's worth noting that ESSs can provide upward reserves by reducing the charging power. Similarly, downward reserve capacity can be offered by charging more from their predetermined power schedules.

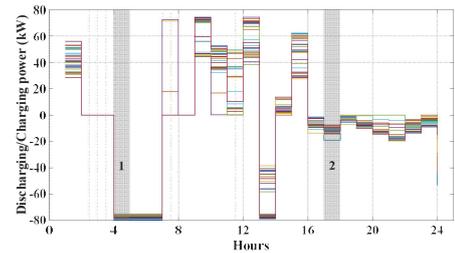

(a) Charging/discharging power of storages
(Discharging corresponds to positive values, while negative values for charging)

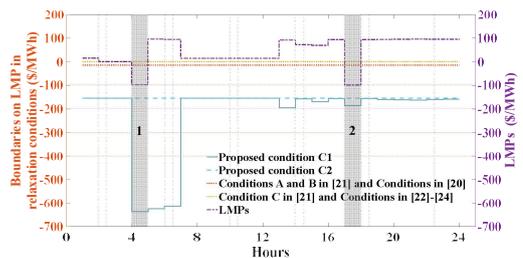

(b) Verifying the validity of relaxation conditions

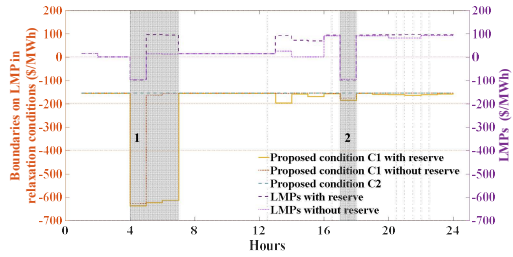

(c) Comparisons of relaxation conditions' boundaries with and without providing reserve services

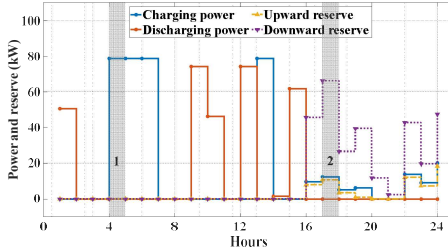

(d) Upward/downward reserves supplied by the storage at bus 10

Fig. 5. The results employing the relaxed model under negative price scenarios.

## C. Verifying the exactness of the proposed sufficient conditions in Theorem 1 under multiple scenarios

In the following scenarios, we assume that the renewable generation cost is -5 $/MWh, which represents the cases of high penetrations of renewable energy or transmission congestion. Then, we compare the results obtained by the proposed relaxation methods in ***Theorem 1*** with those of the MIP model. Table II shows that the relaxed model can obtain consistent energy storage charging/discharging power and objective values as those of the MIP model in all six scenarios, which justifies the validity of the proposed sufficient conditions. The results also show the relaxed model can be calculated much more efficiently than the original MIP model even for a relatively small system. In fact, the effect becomes more significant as both the system scale and the quantity of storage devices increase. This is almost self-evident since the proposed method can eliminate more integer variables in such scenarios.

TABLE II
VERIFYING THE ACCURACY AND EFFICIENCY OF THE PROPOSED SUFFICIENT CONDITIONS IN THEOREM 1 UNDER DIFFERENT SCENARIOS

| Scenario | S1 | S2 | S3 | S4 | S5 | S6 |
| --- | --- | --- | --- | --- | --- | --- |
| Charging price ($/MWh) | 0 | -15 | -10 | 0 | 5 | 10 |
| Discharging price ($/MWh) | 0 | 20 | 15 | 15 | 15 | 15 |
| Whether the Condition C1 is met | yes | yes | yes | yes | yes | yes |
| Whether the Condition C2 is met | yes | yes | yes | yes | yes | yes |
| Whether the relaxed model is accurate | yes | yes | yes | yes | yes | yes |
| Optimal value of the relaxed model ($) | 3776.8 | 3790.4 | 3807.4 | 3941.6 | 3991.3 | 4034.6 |
| Optimal value of the MIP model ($) | 3776.8 | 3790.4 | 3807.4 | 3941.6 | 3991.3 | 4034.6 |
| Solving time of the relaxed model (s) | 14.83 | 12.75 | 12.35 | 14.61 | 11.80 | 11.82 |
| Calculation time of the MIP model (s) | 75.48 | 68.74 | 64.93 | 218.82 | 164.57 | 63.21 |

***Discussion 1***: The cost terms associated with storages in the objective function represent at least three storage dispatch scenarios, involving several combinations of signs (positive or negative) for both discharging prices $g'_{n,t}$ and charging prices $f'_{n,t}$, as detailed in Table III:

TABLE III
THREE APPLICATION SCENARIOS WHERE THE PROPOSED RELAXATION CONDITIONS CAN BE APPLIED

| Scenario | Description | Signs of ($g'_{n,t}, f'_{n,t}$) |
| --- | --- | --- |
| Scenario 1 | The operational costs of ESSs are disregarded in the grid's dispatch at any time $t$ | $g'_{n,t} = 0, f'_{n,t} = 0$ |
| Scenario 2 | The storages' owner pays the grid for the charging energy, and in turn, the grid compensates the ESS owner for the discharging energy at time t | $g'_{n,t} > 0, f'_{n,t} < 0$ |
| Scenario 3 | The grid pays the ESS owner for both charging and discharging energy at time t | $g'_{n,t} > 0, f'_{n,t} > 0$ |

1) In scenario 1, the costs related to charging and discharging are omitted, as discussed in [8], [19], [36], [37]. This situation may happen when the storage device is owned by the power grid company itself. S1 in Table II falls into this category.

2) Scenario 2, which has been extensively considered in the literature [20]-[24], involves the storage owner paying the grid for charging energy while receiving compensation from the grid for discharging. To incentivize owners to discharge stored energy back to the grid, the marginal compensation offered for discharging one unit of energy should cover the owner's marginal charging cost for the equivalent amount of energy, i.e. $g'_{n,t} \geq |f'_{n,t}|$. Scenarios S2 and S3 in Table II are representative examples.

3) Scenario 3 is likely to occur when storages facilitate the absorption of excess renewable energy and help to maintain the grid's energy balance. In this scenario, the grid incentivizes storage charging. Note that such a situation may occur in practice but not necessarily in every period. Scenarios S5 and S6 in Table II fall within this category.

In summary, the storage-concerned OPF model described above is a versatile and generic one, which can be applied to multiple common scenarios. Thus, the derived relaxation conditions based on this general model are generally applicable.

## D. Verifying the validity of different sufficient conditions under negative LMP scenarios

Under different charging/discharging prices and renewable generation costs, we verify whether the proposed sufficient conditions and other sufficient conditions in [8], [20]-[24] are valid. The results are presented in Table IV.

We can find that the proposed relaxation conditions still hold even with a very negative LMP, as shown in scenarios S7-S11. Therefore, we can draw the conclusion that as long as the proposed relaxation conditions hold, it can be guaranteed that there is no SCD. In contrast, those relaxation conditions given in [8], [20]-[24] do not hold in these scenarios, which highlights the advantages of the proposed conditions. Moreover, the numerical tests also validate the correctness of *Lemma 2* and demonstrate that the proposed sufficient conditions have a broader application range compared with the other existing sufficient conditions.

TABLE IV
VERIFYING THE VALIDITY OF DIFFERENT SUFFICIENT CONDITIONS UNDER NEGATIVE LMP SCENARIOS

| Scenario | S7 | S8 | S9 | S10 | S11 |
|---|---|---|---|---|---|
| Charging price ($/MWh) | -8 | -5 | -3 | 15 | 20 |
| Discharging price ($/MWh) | 15 | 15 | 15 | 15 | 15 |
| Renewable generation cost ($/MWh) | -10 | -20 | -30 | -100 | -120 |
| Minimal LMP ($/MWh) | -10.67 | -20.53 | -30.28 | -99.1 | -118.8 |
| Whether the conditions in [8] are met | no | no | no | no | no |
| Whether the conditions in [20] are met | no | no | no | no | no |
| Whether the conditions of group B in [21] are met | no | no | no | no | no |
| Whether the conditions of group C in [21] are met | no | no | no | no | no |
| Whether the conditions in [22] are met | no | no | no | no | no |
| Whether the conditions in [23] are met | no | no | no | no | no |
| Whether the conditions in [24] are met | no | no | no | no | no |
| Whether Condition C1 is met | yes | yes | yes | yes | yes |
| Whether Condition C2 is met | yes | yes | yes | yes | yes |
| Whether the relaxed model is accurate | yes | yes | yes | yes | yes |

During periods when the available renewable energy exceeds the load demand, renewable energy producers tend to offer electricity at negative prices to promote the integration of renewable energy. Nevertheless, in these scenarios, $\underline{\text{LMP}}_{j,t}$ is significantly influenced by the cost coefficient of renewable energy. Considering that the derived relaxation conditions are related to the charging and discharging prices, thus we present the exact relaxation range for storages under varying charging price and renewable energy cost coefficient, as illustrated in Fig. 6. For the sake of illustration, we unify the charging price and the cost coefficient for renewable energy. Moreover, the discharging price is set as 15 $/MWh, with positive generation cost coefficients and reserve cost coefficients given as mentioned. We can observe that:

1) The exact relaxation range obtained based on the proposed *Condition C2* is much larger than the others in [20]-[21]. As for those conditions in [22]-[24], they are all included in the *relaxation conditions of group C in* [21], as demonstrated by the proofs in *Appendix D*. Therefore, we can conclude that the region about LMPs outlined by the proposed *Condition C2* includes those described in [20]-[24], validating the correctness of *Lemma 2*.

2) In this figure, the solid red line and brown line correspond to the cases of $\sigma_{\text{ESS}} = 3$ and $\sigma_{\text{ESS}} = 30$, respectively. It can be found that the exact relaxation boundary gradually shifts downward as $\sigma_{\text{ESS}}$ increases, that is, an increase in the penalty coefficient $\sigma_{\text{ESS}}$ leads to a broader relaxation range. This insight suggests that an increase in the penalty coefficient $\sigma_{\text{ESS}}$, if appropriately adjusted, can guarantee the exactness of the relaxation, validating the effectiveness of the procedure presented in Fig. 4.

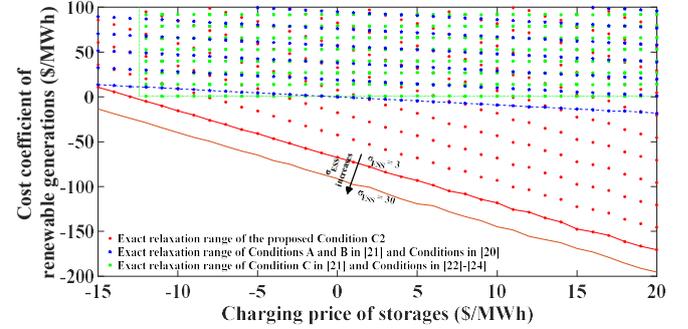

Fig. 6. The exact relaxation range for storages under varying charging price and renewable energy cost coefficient.

## V. CONCLUSION

In this paper, two sufficient conditions for the exact relaxation of complementarity constraints are developed that provably guarantee no SCD in relaxed multi-period OPF. We also prove that the regions on LMPs formed by these two conditions both contain the other existing typical ones. Moreover, the applicable prerequisite of relaxation conditions is generalized from positive LMPs to negative price scenarios. In contrast to just involvement in the energy dispatch, we highlight that offering reserve services can extend the applicability of the proposed relaxation conditions for storages.

The accuracy and advantages of the proposed conditions are demonstrated in numerical tests. It is worth noting that whether the proposed relaxation condition holds can be preliminarily assessed by employing the LMP prediction, facilitating its practical applications. The numerical results also show the broad application prospects of the proposed relaxation methods in the operations of power grid with a high penetration of renewable energy and ESSs.

Moreover, based on the proposed conditions, the duality theory can be employed to perform in-depth theoretical analysis. Correspondingly, it facilitates the optimization of ESSs to be embedded in more complicated problems, for instance, multi-level optimization problems [38] and robust optimization problems [39].

# Supplementary File of Sufficient Conditions for the Exact Relaxation of Complementarity Constraints for Storages in Multi-period OPF Problems

Qi Wang, Wenchuan Wu, *Fellow, IEEE*, Chenhui Lin, Shuwei Xu, Xueliang Li

For the sake of description, we present the mathematical model detailed in the main text:

### A. Objective function

$$\min \ obj = \sum_{t=1}^{T} \left\{ \sum_{n \in \Phi_{\text{ESS}}} \left[ g_{n,t}\left(p_{\text{ESS},n,t}^{\text{dc}}\right) + f_{n,t}\left(p_{\text{ESS},n,t}^{\text{ch}}\right) \right] \right.$$
$$+ \sum_{g \in \Phi_{\text{G}}} C_{g,t}\left(p_{g,t}^{\text{G}}\right) + \sum_{n \in \Phi_{\text{RG}}} C_{n,t}\left(p_{n,t}^{\text{RG}}\right)$$
$$+ \sum_{n \in \Phi_{\text{ESS}}} C_{n,t}^{\text{R}}\left(R_{n,t}^{\text{ch,ESS}-}, R_{n,t}^{\text{dc,ESS}-}, R_{n,t}^{\text{ch,ESS}+}, R_{n,t}^{\text{dc,ESS}+}\right)$$
$$\left. + \sum_{g \in \Phi_{\text{G}}} C_{g,t}^{\text{R}}\left(R_{g,t}^{\text{G}+}, R_{g,t}^{\text{G}-}\right) \right\} + C_{\text{ESS}}^{\text{Pen}} \tag{1}$$

$$C_{g,t}\left(p_{g,t}^{\text{G}}\right) = a_{2,g}\left(p_{g,t}^{\text{G}}\right)^2 + a_{1,g} p_{g,t}^{\text{G}} + a_{0,g} \tag{2}$$

$$C_{n,t}\left(p_{n,t}^{\text{RG}}\right) = b_{n,t} p_{n,t}^{\text{RG}} + \frac{\sigma_{\text{RG}}\left(p_{n,t}^{\text{RG}} - \tilde{P}_{n,t}^{\text{RG}}\right)^2}{\tilde{P}_{n,t}^{\text{RG}}} \tag{3}$$

$$C_{g,t}^{\text{R}}\left(R_{g,t}^{\text{G}+}, R_{g,t}^{\text{G}-}\right) = c_{g,t}^{\text{G}+} \cdot R_{g,t}^{\text{G}+} + c_{g,t}^{\text{G}-} \cdot R_{g,t}^{\text{G}-} \tag{4}$$

$$C_{n,t}^{\text{R}}\left(R_{n,t}^{\text{ch,ESS}-}, R_{n,t}^{\text{dc,ESS}-}, R_{n,t}^{\text{ch,ESS}+}, R_{n,t}^{\text{dc,ESS}+}\right)$$
$$= c_{n,t}^{\text{ESS}+} \cdot \left(R_{n,t}^{\text{ch,ESS}+} + R_{n,t}^{\text{dc,ESS}+}\right) + c_{n,t}^{\text{ESS}-} \cdot \left(R_{n,t}^{\text{ch,ESS}-} + R_{n,t}^{\text{dc,ESS}-}\right) \tag{5}$$

$$C_{\text{ESS}}^{\text{Pen}} = \sigma_{\text{ESS}} \sum_{t=1}^{T} \sum_{n \in \Phi_{\text{ESS}}} \left\{ p_{\text{ESS},n,t}^{\text{dc}}\left(\frac{1}{\eta_{\text{ESS},n}^{\text{dc}}} - 1\right) + p_{\text{ESS},n,t}^{\text{ch}}\left(1 - \eta_{\text{ESS},n}^{\text{ch}}\right) \right\} \tag{6}$$

### B. Operational constraints

$$P_{ij,t} = \frac{1}{\tau_{ij,t}^2} g_{ij}^{\varepsilon} V_{i,t}^2 - \frac{1}{\tau_{ij,t}} V_{i,t} V_{j,t} \left[ g_{ij}^{\varepsilon} \cos\left(\theta_{i,t} - \theta_{j,t} - \phi_{ij,t}\right) \right.$$
$$\left. + b_{ij}^{\varepsilon} \sin\left(\theta_{i,t} - \theta_{j,t} - \phi_{ij,t}\right) \right], \forall ij \in \mathbb{E}, \forall t \tag{7}$$

$$P_{ji,t} = g_{ij}^{\varepsilon} V_{j,t}^2 - \frac{1}{\tau_{ij,t}} V_{i,t} V_{j,t} \left[ g_{ij}^{\varepsilon} \cos\left(\theta_{j,t} - \theta_{i,t} + \phi_{ij,t}\right) \right.$$
$$\left. + b_{ij}^{\varepsilon} \sin\left(\theta_{j,t} - \theta_{i,t} + \phi_{ij,t}\right) \right], \forall ij \in \mathbb{E}, \forall t \tag{8}$$

$$Q_{ij,t} = -\frac{1}{\tau_{ij,t}^2}\left(b_{ij}^{\varepsilon} + \frac{b_{ij}^C}{2}\right) V_{i,t}^2 - \frac{1}{\tau_{ij,t}} V_{i,t} V_{j,t} \left[ g_{ij}^{\varepsilon} \sin\left(\theta_{i,t} - \theta_{j,t} - \phi_{ij,t}\right) \right.$$
$$\left. - b_{ij}^{\varepsilon} \cos\left(\theta_{i,t} - \theta_{j,t} - \phi_{ij,t}\right) \right], \forall ij \in \mathbb{E}, \forall t \tag{9}$$

$$Q_{ji,t} = -\left(b_{ij}^{\varepsilon} + \frac{b_{ij}^C}{2}\right) V_{j,t}^2 - \frac{1}{\tau_{ij,t}} V_{i,t} V_{j,t} \left[ g_{ij}^{\varepsilon} \sin\left(\theta_{j,t} - \theta_{i,t} + \phi_{ij,t}\right) \right.$$
$$\left. - b_{ij}^{\varepsilon} \cos\left(\theta_{j,t} - \theta_{i,t} + \phi_{ij,t}\right) \right], \forall ij \in \mathbb{E}, \forall t \tag{10}$$

$$p_{j,t} - \sum_{i: ji \in \mathbb{E}} P_{ji,t} - \sum_{i: ij \in \mathbb{E}} P_{ji,t} - V_{j,t}^2 g_j^s = 0, \forall j \in \mathcal{N}, \forall t \quad \left(\lambda_{j,t}^{\text{P}}\right) \tag{11}$$

$$q_{j,t} - \sum_{i: ji \in \mathbb{E}} Q_{ji,t} - \sum_{i: ij \in \mathbb{E}} Q_{ji,t} + V_{j,t}^2 b_j^s = 0, \forall j \in \mathcal{N}, \forall t \tag{12}$$

$$p_{j,t} = \sum_{g \in \Phi_{\text{G},j}} p_{g,t}^{\text{G}} + \sum_{n \in \Phi_{\text{RG},j}} p_{n,t}^{\text{RG}} + \sum_{n \in \Phi_{\text{ESS},j}} \left(p_{\text{ESS},n,t}^{\text{dc}} - p_{\text{ESS},n,t}^{\text{ch}}\right) - p_{j,t}^{\text{D}} \tag{13}$$

$$q_{j,t} = \sum_{g \in \Phi_{\text{G},j}} q_{g,t}^{\text{G}} + \sum_{n \in \Phi_{\text{RG},j}} q_{n,t}^{\text{RG}} + \sum_{n \in \Phi_{\text{ESS},j}} q_{\text{ESS},n,t} + \sum_{n \in \Phi_{\text{SVC},j}} q_{n,t}^{\text{SVC}} - q_{j,t}^{\text{D}} \tag{14}$$

$$-p_{\text{ESS},n,t}^{\text{ch}} \le 0, p_{\text{ESS},n,t}^{\text{ch}} \le \overline{P}_{\text{ESS},n}^{\text{ch}}, \quad \left(\lambda_{\text{ESS},n,t}^{\text{ch},1}, \lambda_{\text{ESS},n,t}^{\text{ch},2}\right) \tag{15}$$

$$-p_{\text{ESS},n,t}^{\text{dc}} \le 0, p_{\text{ESS},n,t}^{\text{dc}} \le \overline{P}_{\text{ESS},n}^{\text{dc}}, \quad \left(\lambda_{\text{ESS},n,t}^{\text{dc},1}, \lambda_{\text{ESS},n,t}^{\text{dc},2}\right) \tag{16}$$

$$p_{\text{ESS},n,t}^{\text{ch}} + R_{n,t}^{\text{ch,ESS}-} \le \overline{P}_{\text{ESS},n}^{\text{ch}}, \quad \left(\lambda_{\text{ESS},n,t}^{\text{RD},1}\right)$$
$$R_{n,t}^{\text{dc,ESS}-} \le p_{\text{ESS},n,t}^{\text{dc}}, \quad \left(\lambda_{\text{ESS},n,t}^{\text{RD},2}\right) \tag{17}$$

$$p_{\text{ESS},n,t}^{\text{dc}} + R_{n,t}^{\text{dc,ESS}+} \le \overline{P}_{\text{ESS},n}^{\text{dc}}, \quad \left(\lambda_{\text{ESS},n,t}^{\text{RU},1}\right)$$
$$R_{n,t}^{\text{ch,ESS}+} \le p_{\text{ESS},n,t}^{\text{ch}}, \quad \left(\lambda_{\text{ESS},n,t}^{\text{RU},2}\right) \tag{18}$$

$$R_{n,t}^{\text{ch,ESS}-}, R_{n,t}^{\text{dc,ESS}-}, R_{n,t}^{\text{ch,ESS}+}, R_{n,t}^{\text{dc,ESS}+} \ge 0 \tag{19}$$

$$\left(1-\delta_{\text{ESS},n}\right)^t E_{\text{ESS},n,0} + \sum_{\tau=1}^t \left(1-\delta_{\text{ESS},n}\right)^{t-\tau}$$
$$\times \left(\left(p_{\text{ESS},n,\tau}^{\text{ch}} - R_{n,t}^{\text{ch,ESS}-}\right) \cdot \eta_{\text{ESS},n}^{\text{ch}} - \left(p_{\text{ESS},n,\tau}^{\text{dc}} + R_{n,t}^{\text{dc,ESS}+}\right)/\eta_{\text{ESS},n}^{\text{dc}}\right) \Delta t \ge \underline{E}_{\text{ESS},n}, \left(\lambda_{\text{ESS},n,t}^{\text{SOC},1}\right) \tag{20}$$

$$\left(1-\delta_{\text{ESS},n}\right)^t E_{\text{ESS},n,0} + \sum_{\tau=1}^t \left(1-\delta_{\text{ESS},n}\right)^{t-\tau}$$
$$\times \left(\left(p_{\text{ESS},n,\tau}^{\text{ch}} + R_{n,t}^{\text{ch,ESS}-}\right) \cdot \eta_{\text{ESS},n}^{\text{ch}} - \left(p_{\text{ESS},n,\tau}^{\text{dc}} - R_{n,t}^{\text{dc,ESS}-}\right)/\eta_{\text{ESS},n}^{\text{dc}}\right) \Delta t \le \overline{E}_{\text{ESS},n}, \left(\lambda_{\text{ESS},n,t}^{\text{SOC},2}\right) \tag{21}$$

$$E_{\text{ESS},n,T} = \left(1-\delta_{\text{ESS},n}\right)^T E_{\text{ESS},n,0} + \sum_{\tau=1}^T \left(1-\delta_{\text{ESS},n}\right)^{T-\tau}$$
$$\times \left(p_{\text{ESS},n,\tau}^{\text{ch}} \cdot \eta_{\text{ESS},n}^{\text{ch}} - p_{\text{ESS},n,\tau}^{\text{dc}}/\eta_{\text{ESS},n}^{\text{dc}}\right) \Delta t = E_{\text{ESS},n,0} \tag{22}$$

$$\left(p_{\text{ESS},n,t}^{\text{dc}}\right)^2 + \left(q_{\text{ESS},n,t}\right)^2 \le \left(S_{\text{ESS},n}\right)^2, \quad \left(\lambda_{\text{ESS},n,t}^{\text{S},1}\right) \tag{23}$$

$$\left(p_{\text{ESS},n,t}^{\text{ch}}\right)^2 + \left(q_{\text{ESS},n,t}\right)^2 \le \left(S_{\text{ESS},n}\right)^2, \quad \left(\lambda_{\text{ESS},n,t}^{\text{S},2}\right) \tag{24}$$

$$\frac{p_{\text{ESS},n,t}^{\text{ch}}}{\overline{P}_{\text{ESS},n}^{\text{ch}}} + \frac{p_{\text{ESS},n,t}^{\text{dc}}}{\overline{P}_{\text{ESS},n}^{\text{dc}}} \le 1, \forall n \in \Phi_{\text{ESS}}, \forall t, \quad \left(\lambda_{\text{ESS},n,t}^{\text{relax}}\right) \tag{25}$$

$$\underline{P}_g^{\text{G}} \le p_{g,t}^{\text{G}} \le \overline{P}_g^{\text{G}}, \forall g \in \Phi_{\text{G}} \tag{26}$$

$$\underline{Q}_g^{\text{G}} \le q_{g,t}^{\text{G}} \le \overline{Q}_g^{\text{G}}, \forall g \in \Phi_{\text{G}} \tag{27}$$

$$-RD_g \Delta t \leq p_{g,t+1}^{G} - p_{g,t}^{G} \leq RU_g \Delta t \tag{28}$$

$$0 \leq R_{g,t}^{G+} \leq RU_g \Delta t, \tag{29}$$
$$R_{g,t}^{G+} \leq \overline{P}_g^{G} - p_{g,t}^{G}, \forall g \in \Phi_G$$

$$0 \leq R_{g,t}^{G-} \leq RD_g \Delta t, \tag{30}$$
$$R_{g,t}^{G-} \leq p_{g,t}^{G} - \underline{P}_g^{G}, \forall g \in \Phi_G$$

$$\sum_{g \in \Phi_G} R_{g,t}^{G+} + \sum_{n \in \Phi_{ESS}} \left( R_{n,t}^{ch,ESS+} + R_{n,t}^{dc,ESS+} \right) \geq SRU_t, \tag{31}$$
$$\sum_{g \in \Phi_G} R_{g,t}^{G-} + \sum_{n \in \Phi_{ESS}} \left( R_{n,t}^{ch,ESS-} + R_{n,t}^{dc,ESS-} \right) \geq SRD_t$$

$$\underline{P}_{n,t}^{RG} \leq p_{n,t}^{RG} \leq \widetilde{P}_{n,t}^{RG}, \forall n \in \Phi_{RG} \tag{32}$$

$$\left(p_{n,t}^{RG}\right)^2 + \left(q_{n,t}^{RG}\right)^2 \leq \left(S_n^{RG}\right)^2, \forall n \in \Phi_{RG} \tag{33}$$

$$\underline{Q}_n^{SVC} \leq q_{n,t}^{SVC} \leq \overline{Q}_n^{SVC}, \forall n \in \Phi_{SVC} \tag{34}$$

$$\underline{V}_i \leq V_{i,t} \leq \overline{V}_i, \forall i \in \mathcal{N} \tag{35}$$

$$P_{ij,t}^2 + Q_{ij,t}^2 \leq \overline{S}_{ij}^2, \forall ij \in \mathbb{E} \tag{36}$$

## APPENDIX A
### PROOF FOR THEOREM 1

**Proof of Condition C1**: Denoting the Lagrangian function of the relaxed model as $L_{RM}$, and defining a new variable $\Gamma_{ESS,n,t} = \sum_{\tau=t}^{T} \left(1 - \delta_{ESS,n}\right)^{\tau-t} \left( \lambda_{ESS,n,\tau}^{SOC,2} - \lambda_{ESS,n,\tau}^{SOC,1} \right)$. The constraints related to $p_{ESS,n,t}^{ch}$ and $p_{ESS,n,t}^{dc}$ include (11), (13)-(18), (20)-(25). Then employing KKT conditions, the following equations are obtained:

$$\frac{\partial L_{RM}}{\partial \left(p_{ESS,n,t}^{ch}\right)} = \frac{\partial obj}{\partial \left(p_{ESS,n,t}^{ch}\right)} + \lambda_{j,t}^{p} - \lambda_{ESS,n,t}^{ch,1} + \lambda_{ESS,n,t}^{ch,2} + \frac{\lambda_{ESS,n,t}^{relax}}{\overline{P}_{ESS,n}^{ch}} \tag{37}$$
$$+ 2\lambda_{ESS,n,t}^{S,2} \cdot p_{ESS,n,t}^{ch} + \eta_{ESS,n}^{ch} \cdot \Gamma_{ESS,n,t} \cdot \Delta t + \lambda_{ESS,n,t}^{RD,1} = 0$$

$$\frac{\partial L_{RM}}{\partial \left(p_{ESS,n,t}^{dc}\right)} = \frac{\partial obj}{\partial \left(p_{ESS,n,t}^{dc}\right)} - \lambda_{j,t}^{p} - \lambda_{ESS,n,t}^{dc,1} + \lambda_{ESS,n,t}^{dc,2} + \frac{\lambda_{ESS,n,t}^{relax}}{\overline{P}_{ESS,n}^{dc}} \tag{38}$$
$$+ 2\lambda_{ESS,n,t}^{S,1} \cdot p_{ESS,n,t}^{dc} - \Gamma_{ESS,n,t} \cdot \Delta t / \eta_{ESS,n}^{dc} + \lambda_{ESS,n,t}^{RU,1} = 0$$

With $(37) / \eta_{ESS,n}^{dc} + (38) \cdot \eta_{ESS,n}^{ch}$, we have

$$\left[ \frac{\partial obj}{\partial \left(p_{ESS,n,t}^{ch}\right)} + \lambda_{j,t}^{p} - \lambda_{ESS,n,t}^{ch,1} + \lambda_{ESS,n,t}^{ch,2} + 2\lambda_{ESS,n,t}^{S,2} \cdot p_{ESS,n,t}^{ch} \right.$$
$$\left. + \frac{\lambda_{ESS,n,t}^{relax}}{\overline{P}_{ESS,n}^{ch}} + \lambda_{ESS,n,t}^{RD,1} \right] / \eta_{ESS,n}^{dc} + \left[ \frac{\partial obj}{\partial \left(p_{ESS,n,t}^{dc}\right)} - \lambda_{j,t}^{p} - \lambda_{ESS,n,t}^{dc,1} \right. \tag{39}$$
$$\left. + \lambda_{ESS,n,t}^{dc,2} + 2\lambda_{ESS,n,t}^{S,1} \cdot p_{ESS,n,t}^{dc} + \frac{\lambda_{ESS,n,t}^{relax}}{\overline{P}_{ESS,n}^{dc}} + \lambda_{ESS,n,t}^{RU,1} \right] \cdot \eta_{ESS,n}^{ch} = 0$$

Further, equation (39) yields:

$$\lambda_{j,t}^{p} = \left[ \frac{\partial obj}{\partial \left(p_{ESS,n,t}^{ch}\right)} / \eta_{ESS,n}^{dc} + \frac{\partial obj}{\partial \left(p_{ESS,n,t}^{dc}\right)} \cdot \eta_{ESS,n}^{ch} \right] / \left( \eta_{ESS,n}^{ch} - \frac{1}{\eta_{ESS,n}^{dc}} \right)$$
$$- \left\{ \left[ \lambda_{ESS,n,t}^{ch,2} + 2\lambda_{ESS,n,t}^{S,2} \cdot p_{ESS,n,t}^{ch} + \frac{\lambda_{ESS,n,t}^{relax}}{\overline{P}_{ESS,n}^{ch}} + \lambda_{ESS,n,t}^{RD,1} \right] / \eta_{ESS,n}^{dc} \right.$$
$$\left. + \left[ \lambda_{ESS,n,t}^{dc,2} + 2\lambda_{ESS,n,t}^{S,1} \cdot p_{ESS,n,t}^{dc} + \frac{\lambda_{ESS,n,t}^{relax}}{\overline{P}_{ESS,n}^{dc}} + \lambda_{ESS,n,t}^{RU,1} \right] \cdot \eta_{ESS,n}^{ch} \right\} / \left( \frac{1}{\eta_{ESS,n}^{dc}} - \eta_{ESS,n}^{ch} \right)$$
$$+ \left( \lambda_{ESS,n,t}^{ch,1} / \eta_{ESS,n}^{dc} + \lambda_{ESS,n,t}^{dc,1} \cdot \eta_{ESS,n}^{ch} \right) / \left( 1/\eta_{ESS,n}^{dc} - \eta_{ESS,n}^{ch} \right) \tag{40}$$

If the proposed *Condition C1* holds, formula (41) can be obtained from (40).

$$\left( \lambda_{ESS,n,t}^{ch,1} / \eta_{ESS,n}^{dc} + \lambda_{ESS,n,t}^{dc,1} \cdot \eta_{ESS,n}^{ch} \right) / \left( 1/\eta_{ESS,n}^{dc} - \eta_{ESS,n}^{ch} \right) > 0 \tag{41}$$

Since $0 < \eta_{ESS,n}^{dc}, \eta_{ESS,n}^{ch} < 1$, we have

$$\lambda_{ESS,n,t}^{ch,1} / \eta_{ESS,n}^{dc} + \lambda_{ESS,n,t}^{dc,1} \cdot \eta_{ESS,n}^{ch} > 0 \tag{42}$$

Based on (42), we can further derive inequality (43) in consideration of $\eta_{ESS,n}^{dc}, \eta_{ESS,n}^{ch} > 0$ and $\lambda_{ESS,n,t}^{ch,1}, \lambda_{ESS,n,t}^{dc,1} \geq 0$.

$$\lambda_{ESS,n,t}^{ch,1} + \lambda_{ESS,n,t}^{dc,1} > 0 \tag{43}$$

(43) indicates that at least one of $p_{ESS,n,t}^{ch} = 0$ and $p_{ESS,n,t}^{dc} = 0$ holds considering the complementary slackness conditions of constraints (15) and (16), which means that $p_{ESS,n,t}^{ch} \cdot p_{ESS,n,t}^{dc} = 0$. Hence, under this condition, even if the complementary constraints are relaxed, the optimal solution still implicitly avoids the SCD.
∎

**Proof of Condition C2**: In combination with the expression of $\lambda_{j,t}^{p}$ in equation (40), we can conclude that if the proposed *Condition C2* holds, we can obtain (44).

$$\left( \lambda_{ESS,n,t}^{ch,1} / \eta_{ESS,n}^{dc} + \lambda_{ESS,n,t}^{dc,1} \cdot \eta_{ESS,n}^{ch} \right) / \left( 1/\eta_{ESS,n}^{dc} - \eta_{ESS,n}^{ch} \right) >$$
$$\left\{ \left[ \lambda_{ESS,n,t}^{ch,2} + 2\lambda_{ESS,n,t}^{S,2} \cdot p_{ESS,n,t}^{ch} + \frac{\lambda_{ESS,n,t}^{relax}}{\overline{P}_{ESS,n}^{ch}} + \lambda_{ESS,n,t}^{RD,1} \right] / \eta_{ESS,n}^{dc} \right.$$
$$\left. + \left[ \lambda_{ESS,n,t}^{dc,2} + 2\lambda_{ESS,n,t}^{S,1} \cdot p_{ESS,n,t}^{dc} + \frac{\lambda_{ESS,n,t}^{relax}}{\overline{P}_{ESS,n}^{dc}} + \lambda_{ESS,n,t}^{RU,1} \right] \cdot \eta_{ESS,n}^{ch} \right\} / \left( 1/\eta_{ESS,n}^{dc} - \eta_{ESS,n}^{ch} \right) \tag{44}$$

Considering $\eta_{ESS,n}^{dc}, \eta_{ESS,n}^{ch} > 0$, $\lambda_{ESS,n,t}^{ch,2}, \lambda_{ESS,n,t}^{dc,2} \geq 0$ and $\lambda_{ESS,n,t}^{S,1}, \lambda_{ESS,n,t}^{S,2}, \lambda_{ESS,n,t}^{relax}, \lambda_{ESS,n,t}^{RD,1}, \lambda_{ESS,n,t}^{RU,1}, p_{ESS,n,t}^{ch}, p_{ESS,n,t}^{dc} \geq 0$, we can derive

$$\left[ \lambda_{ESS,n,t}^{ch,2} + 2\lambda_{ESS,n,t}^{S,2} \cdot p_{ESS,n,t}^{ch} + \frac{\lambda_{ESS,n,t}^{relax}}{\overline{P}_{ESS,n}^{ch}} + \lambda_{ESS,n,t}^{RD,1} \right] / \eta_{ESS,n}^{dc}$$
$$+ \left[ \lambda_{ESS,n,t}^{dc,2} + 2\lambda_{ESS,n,t}^{S,1} \cdot p_{ESS,n,t}^{dc} + \frac{\lambda_{ESS,n,t}^{relax}}{\overline{P}_{ESS,n}^{dc}} + \lambda_{ESS,n,t}^{RU,1} \right] \cdot \eta_{ESS,n}^{ch} \geq 0 \tag{45}$$

Combining formulas (44), (45) and $1/\eta_{ESS,n}^{dc} > \eta_{ESS,n}^{ch}$, we have

$$\left( \lambda_{ESS,n,t}^{ch,1} / \eta_{ESS,n}^{dc} + \lambda_{ESS,n,t}^{dc,1} \cdot \eta_{ESS,n}^{ch} \right) > 0 \tag{46}$$

Since $\eta_{ESS,n}^{dc}, \eta_{ESS,n}^{ch} > 0$ and $\lambda_{ESS,n,t}^{ch,1}, \lambda_{ESS,n,t}^{dc,1} \geq 0$, (47) holds.

$$\lambda_{ESS,n,t}^{ch,1} + \lambda_{ESS,n,t}^{dc,1} > 0 \tag{47}$$

Therefore, *Condition C2* indicates that $p_{ESS,n,t}^{ch} \cdot p_{ESS,n,t}^{dc} = 0$. In other words, based on *Condition C2*, the optimal solution enforces the complementary constraint, even though it is relaxed.
∎

## APPENDIX B
## PROOF FOR LEMMA 1

*Proof*: Combining formulas (37) and (38), then **Lemma 1** can be equivalently described as follows:

$$\lambda_{j,t}^{p,C2} \geq \lambda_{j,t}^{p,C1} \qquad (48)$$

Further, combined with formulas (39) and (40), we can find that formula (48) holds if and only if formula (49) holds.

$$\left\{ \left[ \lambda_{ESS,n,t}^{ch,2} + 2\lambda_{ESS,n,t}^{S,2} \cdot p_{ESS,n,t}^{ch} + \frac{\lambda_{ESS,n,t}^{relax}}{\overline{P}_{ESS,n}^{ch}} + \lambda_{ESS,n,t}^{RD,1} \right] / \eta_{ESS,n}^{dc} \right.$$
$$\left. + \left[ \lambda_{ESS,n,t}^{dc,2} + 2\lambda_{ESS,n,t}^{S,1} \cdot p_{ESS,n,t}^{dc} + \frac{\lambda_{ESS,n,t}^{relax}}{\overline{P}_{ESS,n}^{dc}} + \lambda_{ESS,n,t}^{RU,1} \right] \cdot \eta_{ESS,n}^{ch} \right\} / \left( \frac{1}{\eta_{ESS,n}^{dc}} - \eta_{ESS,n}^{ch} \right) \geq 0 \qquad (49)$$

Combining (45) and $1/\eta_{ESS,n}^{dc} > \eta_{ESS,n}^{ch}$, it can be inferred that (49) holds. Thus, we can draw the conclusion that **Lemma 1** holds.

And as long as any one of these constraints corresponding to multipliers $\lambda_{ESS,n,t}^{ch,2}, \lambda_{ESS,n,t}^{dc,2}, \lambda_{ESS,n,t}^{RD,1}, \lambda_{ESS,n,t}^{RU,1}, \lambda_{ESS,n,t}^{S,1}, \lambda_{ESS,n,t}^{S,2}, \lambda_{ESS,n,t}^{relax}$, i.e., constraints (15)-(18) and (23)-(25), becomes active at time $t$, the less than or equal sign in formula (48) will become a strict less than sign.

∎

## APPENDIX C

### C. Proof for the Relaxation Conditions of Group B in [2] under the ACOPF scenario

The **Relaxation conditions of group B** in [2] are given as following. Although they are obtained using the ED model, similar relaxation conditions can be derived in ACOPF scenario.

**C3_1**: $\quad \dfrac{\partial \text{obj}}{\partial \left( p_{ESS,n,t}^{ch} \right)} + \dfrac{\partial \text{obj}}{\partial \left( p_{ESS,n,t}^{dc} \right)} > 0 \qquad (50)$

**C3_2**: $\quad \lambda_{j,t}^{p} \geq \lambda_{j,t}^{p,C3} = -\dfrac{\partial \text{obj}}{\partial \left( p_{ESS,n,t}^{ch} \right)} \qquad (51)$

*Proof*: We prove this type of conditions by contradiction. Suppose there exists $p_{ESS,n,t}^{ch} > 0$ and $p_{ESS,n,t}^{dc} > 0$ for the $n^{th}$ ESS at time $t$ in the optimal solution of the relaxed model. Then we have $\lambda_{ESS,n,t}^{ch,1} = 0$ and $\lambda_{ESS,n,t}^{dc,1} = 0$.

i) Combined with $\lambda_{ESS,n,t}^{ch,1} = 0$ and $\lambda_{ESS,n,t}^{ch,2} \geq 0$, formula (37) yields:

$$\frac{\partial \text{obj}}{\partial \left( p_{ESS,n,t}^{ch} \right)} + \lambda_{j,t}^{p} + \frac{\lambda_{ESS,n,t}^{relax}}{\overline{P}_{ESS,n}^{ch}} + 2\lambda_{ESS,n,t}^{S,2} \cdot p_{ESS,n,t}^{ch} \qquad (52)$$
$$+ \eta_{ESS,n}^{ch} \cdot \Gamma_{ESS,n,t} \cdot \Delta t + \lambda_{ESS,n,t}^{RD,1} \leq 0$$

Furthermore, with $\lambda_{ESS,n,t}^{relax}, \lambda_{ESS,n,t}^{S,2}, \lambda_{ESS,n,t}^{RD,1} \geq 0$, $\eta_{ESS,n}^{ch} > 0$ and **condition C3_2**, we can infer that $\Gamma_{ESS,n,t} \leq 0$.

ii) Summing formulas (37) and (38), and combining $\lambda_{ESS,n,t}^{ch,1} = \lambda_{ESS,n,t}^{dc,1} = 0$, we can obtain

$$\frac{\partial \text{obj}}{\partial \left( p_{ESS,n,t}^{ch} \right)} + \frac{\partial \text{obj}}{\partial \left( p_{ESS,n,t}^{dc} \right)} + \frac{\lambda_{ESS,n,t}^{relax}}{\overline{P}_{ESS,n}^{ch}} + \frac{\lambda_{ESS,n,t}^{relax}}{\overline{P}_{ESS,n}^{dc}} + \lambda_{ESS,n,t}^{ch,2} + \lambda_{ESS,n,t}^{dc,2}$$
$$+ \lambda_{ESS,n,t}^{RD,1} + \lambda_{ESS,n,t}^{RU,1} + 2\left( \lambda_{ESS,n,t}^{S,2} \cdot p_{ESS,n,t}^{ch} + \lambda_{ESS,n,t}^{S,1} \cdot p_{ESS,n,t}^{dc} \right) \qquad (53)$$
$$- \left( 1/\eta_{ESS,n}^{dc} - \eta_{ESS,n}^{ch} \right) \cdot \Gamma_{ESS,n,t} \cdot \Delta t = 0$$

In consideration of **condition C3_1**, $\lambda_{ESS,n,t}^{ch,2}, \lambda_{ESS,n,t}^{dc,2} \geq 0$, $\lambda_{ESS,n,t}^{relax}, \lambda_{ESS,n,t}^{RD,1}, \lambda_{ESS,n,t}^{RU,1}, \lambda_{ESS,n,t}^{S,1}, \lambda_{ESS,n,t}^{S,2} \geq 0$, and $\left( 1/\eta_{ESS,n}^{dc} - \eta_{ESS,n}^{ch} \right) > 0$, it can be concluded that $\Gamma_{ESS,n,t} > 0$, which conflicts with the previous conclusion $\Gamma_{ESS,n,t} \leq 0$.

Therefore, we can draw the conclusion that $p_{ESS,n,t}^{ch} > 0$ and $p_{ESS,n,t}^{dc} > 0$ are impossible to appear in the optimal solution of the relaxed model at the same time for any ESS at any time period on the premise that **condition C3** holds.

∎

### D. Proof for the Relaxation Conditions in [1] and Relaxation Conditions of Group A in [2] under the ACOPF scenario

Similarly, the **Relaxation conditions** in [1] and **Relaxation conditions of group A** in [2] also hold for ACOPF problems. The difference between **condition C4** and **condition C3** lies in their different signs.

**C4_1**: $\quad \dfrac{\partial \text{obj}}{\partial \left( p_{ESS,n,t}^{ch} \right)} + \dfrac{\partial \text{obj}}{\partial \left( p_{ESS,n,t}^{dc} \right)} \geq 0 \qquad (54)$

**C4_2**: $\quad \lambda_{j,t}^{p} > \lambda_{j,t}^{p,C4} = -\dfrac{\partial \text{obj}}{\partial \left( p_{ESS,n,t}^{ch} \right)} \qquad (55)$

*Proof*: We prove this type of conditions by contradiction.

i) With formula (52), $\lambda_{ESS,n,t}^{relax}, \lambda_{ESS,n,t}^{S,2}, \lambda_{ESS,n,t}^{RD,1} \geq 0$, $\eta_{ESS,n}^{ch} > 0$ and **condition C4_2**, we can infer that $\Gamma_{ESS,n,t} < 0$.

ii) Combining formula (53), **condition C4_1**, $\lambda_{ESS,n,t}^{ch,2}, \lambda_{ESS,n,t}^{dc,2} \geq 0$, $\lambda_{ESS,n,t}^{relax}, \lambda_{ESS,n,t}^{RD,1}, \lambda_{ESS,n,t}^{RU,1}, \lambda_{ESS,n,t}^{S,1}, \lambda_{ESS,n,t}^{S,2} \geq 0$ and $\left( 1/\eta_{ESS,n}^{dc} - \eta_{ESS,n}^{ch} \right) > 0$, we can know that $\Gamma_{ESS,n,t} \geq 0$, which is in conflict with the above conclusion $\Gamma_{ESS,n,t} < 0$.

Thus, the relaxation has been proven to be exact for the given **condition C4**.

∎

### E. Proof for the Relaxation Conditions of Group C in [2] under the ACOPF scenario

Under the **Relaxation conditions of Group C** in [2], the ESSs' complementarity constraints are enforced in ED. Here, we give the proof of these conditions in ACOPF scenario:

**C5_1**: $\quad \dfrac{\partial \text{obj}}{\partial \left( p_{ESS,n,t}^{ch} \right)} / \eta_{ESS,n}^{dc} + \dfrac{\partial \text{obj}}{\partial \left( p_{ESS,n,t}^{dc} \right)} \cdot \eta_{ESS,n}^{ch} > 0 \qquad (56)$

**C5_2**: $\quad \lambda_{j,t}^{p} \geq 0 \qquad (57)$

*Proof*: We prove this type of conditions by contradiction. Assume $p_{ESS,n,t}^{ch}, p_{ESS,n,t}^{dc} > 0$, then $\lambda_{ESS,n,t}^{ch,1} = \lambda_{ESS,n,t}^{dc,1} = 0$. With (37) +(38)$\cdot \eta_{ESS,n}^{ch} \cdot \eta_{ESS,n}^{dc}$, we have

$$\left[\frac{\partial \text{obj}}{\partial\left(p_{\text{ESS},n,t}^{\text{ch}}\right)}+\lambda_{j,t}^{\text{p}}-\lambda_{\text{ESS},n,t}^{\text{ch},1}+\lambda_{\text{ESS},n,t}^{\text{ch},2}+2\lambda_{\text{ESS},n,t}^{\text{S},2}\cdot p_{\text{ESS},n,t}^{\text{ch}}\right.$$
$$\left.+\frac{\lambda_{\text{ESS},n,t}^{\text{relax}}}{\overline{P}_{\text{ESS},n}^{\text{ch}}}+\lambda_{\text{ESS},n,t}^{\text{RD},1}\right]+\left[\frac{\partial \text{obj}}{\partial\left(p_{\text{ESS},n,t}^{\text{dc}}\right)}-\lambda_{j,t}^{\text{p}}-\lambda_{\text{ESS},n,t}^{\text{dc},1}+\lambda_{\text{ESS},n,t}^{\text{dc},2}\right. \quad (58)$$
$$\left.+2\lambda_{\text{ESS},n,t}^{\text{S},1}\cdot p_{\text{ESS},n,t}^{\text{dc}}+\frac{\lambda_{\text{ESS},n,t}^{\text{relax}}}{\overline{P}_{\text{ESS},n}^{\text{dc}}}+\lambda_{\text{ESS},n,t}^{\text{RU},1}\right]\cdot \eta_{\text{ESS},n}^{\text{ch}}\cdot \eta_{\text{ESS},n}^{\text{dc}}=0$$

Further, equation (58) yields:

$$\lambda_{j,t}^{\text{p}}\left(1-\eta_{\text{ESS},n}^{\text{ch}}\cdot \eta_{\text{ESS},n}^{\text{dc}}\right)+\frac{\partial \text{obj}}{\partial\left(p_{\text{ESS},n,t}^{\text{ch}}\right)}+\frac{\partial \text{obj}}{\partial\left(p_{\text{ESS},n,t}^{\text{dc}}\right)}\cdot \eta_{\text{ESS},n}^{\text{ch}}\cdot \eta_{\text{ESS},n}^{\text{dc}}$$
$$+\lambda_{\text{ESS},n,t}^{\text{ch},2}+2\lambda_{\text{ESS},n,t}^{\text{S},2}\cdot p_{\text{ESS},n,t}^{\text{ch}}+\frac{\lambda_{\text{ESS},n,t}^{\text{relax}}}{\overline{P}_{\text{ESS},n}^{\text{ch}}}+\lambda_{\text{ESS},n,t}^{\text{RD},1} \quad (59)$$
$$+\left(\lambda_{\text{ESS},n,t}^{\text{dc},2}+2\lambda_{\text{ESS},n,t}^{\text{S},1}\cdot p_{\text{ESS},n,t}^{\text{dc}}+\frac{\lambda_{\text{ESS},n,t}^{\text{relax}}}{\overline{P}_{\text{ESS},n}^{\text{dc}}}+\lambda_{\text{ESS},n,t}^{\text{RU},1}\right)\cdot \eta_{\text{ESS},n}^{\text{ch}}\cdot \eta_{\text{ESS},n}^{\text{dc}}=0$$

Given *condition C5* and $\eta_{\text{ESS},n}^{\text{ch}}\cdot \eta_{\text{ESS},n}^{\text{dc}}<1$, it follows that

$$\lambda_{j,t}^{\text{p}}\left(1-\eta_{\text{ESS},n}^{\text{ch}}\cdot \eta_{\text{ESS},n}^{\text{dc}}\right)+\frac{\partial \text{obj}}{\partial\left(p_{\text{ESS},n,t}^{\text{ch}}\right)}+\frac{\partial \text{obj}}{\partial\left(p_{\text{ESS},n,t}^{\text{dc}}\right)}\cdot \eta_{\text{ESS},n}^{\text{ch}}\cdot \eta_{\text{ESS},n}^{\text{dc}}>0 \quad (60)$$

Further, combining formulas (59) and (60), we have

$$\lambda_{\text{ESS},n,t}^{\text{ch},2}+2\lambda_{\text{ESS},n,t}^{\text{S},2}\cdot p_{\text{ESS},n,t}^{\text{ch}}+\frac{\lambda_{\text{ESS},n,t}^{\text{relax}}}{\overline{P}_{\text{ESS},n}^{\text{ch}}}+\lambda_{\text{ESS},n,t}^{\text{RD},1}$$
$$+\left(\lambda_{\text{ESS},n,t}^{\text{dc},2}+2\lambda_{\text{ESS},n,t}^{\text{S},1}\cdot p_{\text{ESS},n,t}^{\text{dc}}+\frac{\lambda_{\text{ESS},n,t}^{\text{relax}}}{\overline{P}_{\text{ESS},n}^{\text{dc}}}+\lambda_{\text{ESS},n,t}^{\text{RU},1}\right)\cdot \eta_{\text{ESS},n}^{\text{ch}}\cdot \eta_{\text{ESS},n}^{\text{dc}}<0 \quad (61)$$

While as $\lambda_{\text{ESS},n,t}^{\text{ch},2},\lambda_{\text{ESS},n,t}^{\text{dc},2},\lambda_{\text{ESS},n,t}^{\text{S},1},\lambda_{\text{ESS},n,t}^{\text{S},2},\lambda_{\text{ESS},n,t}^{\text{relax}},\lambda_{\text{ESS},n,t}^{\text{RD},1},\lambda_{\text{ESS},n,t}^{\text{RU},1}\geq 0$, there is a contradiction with formula (61). As a conclusion, under the given *condition C5*, the assumption $p_{\text{ESS},n,t}^{\text{ch}},p_{\text{ESS},n,t}^{\text{dc}}>0$ does not hold. Therefore, the relaxation conditions are exact.

∎

*F. Proof for the Relaxation Conditions in [3] under the ACOPF scenario*

The *Relaxation conditions* in [3] are presented as follows:

*C6_1*: $\frac{\partial \text{obj}}{\partial\left(p_{\text{ESS},n,t}^{\text{ch}}\right)}\geq 0,\frac{\partial \text{obj}}{\partial\left(p_{\text{ESS},n,t}^{\text{dc}}\right)}\geq 0,\frac{\partial \text{obj}}{\partial\left(p_{\text{ESS},n,t}^{\text{ch}}\right)}+\frac{\partial \text{obj}}{\partial\left(p_{\text{ESS},n,t}^{\text{dc}}\right)}>0$ (62)

*C6_2*: $\lambda_{j,t}^{\text{p}}\geq 0$ (63)

*Proof*: Comparing *Condition C5* and *Condition C6*, we can find that they differ only in condition *C5_1* and condition *C6_1*. What's more, we can easily derive condition *C5_1* from condition *C6_1*. Thus, we can conclude that the region formulated by *Condition C5* contains that of *Condition C6*. Then, combined with the proof in the previous subsection C in this appendix, we can draw the conclusion that the relaxation is exact for the given *Condition C6* under the ACOPF scenario.

∎

APPENDIX D

PROOF FOR LEMMA 2

In combination of *Lemma 1*, we can infer that if *Lemma 3* holds, we can also conclude that *Lemma 2* is valid.

*Lemma 3*: *The region about LMPs formulated by Condition C2 in Theorem 1 contains those of the other six existing relaxation conditions*:

i) *the Relaxation conditions in* [1] (*i.e., the Relaxation conditions of group A in* [2]);
ii) *the Relaxation conditions of group B in* [2];
iii) *the Relaxation conditions of group C in* [2];
iv) *the Relaxation conditions in* [3];
v) *the Relaxation conditions in* [4];
vi) *the Relaxation conditions in* [5].

In addition to the relaxation conditions in [1]-[3] already presented in Appendix C, we also provide the relaxation conditions in [4]-[5] here to facilitate subsequent proofs.

In [5], the storage penalty term similar to (6) is added in the objective function to avoid SCD. And the following conditions are given:

*C7_1*: $\lambda_{j,t}^{\text{p}}\geq 0$ (64)

*C7_2*: $\lambda_{\text{ESS},n,t}^{\text{S},1}=\lambda_{\text{ESS},n,t}^{\text{S},2}=0$ (65)

The relaxation conditions in [4] are presented as:

*C8_1*: $\frac{\partial \text{obj}}{\partial\left(p_{\text{ESS},n,t}^{\text{ch}}\right)}=0,\frac{\partial \text{obj}}{\partial\left(p_{\text{ESS},n,t}^{\text{dc}}\right)}>0$ (66)

*C8_2*: $\lambda_{j,t}^{\text{p}}>0$ (67)

Subsequently, we redescribe *Lemma 3* as *Lemmas 4-9* equivalently. Correspondingly, if *Lemmas 4-9* all hold, it can be concluded that *Lemma 2* is valid. Below, we will provide the proofs for *Lemmas 4-9*.

*A. Proof for Lemma 4*

*Lemma 4*: *The region about LMPs formulated by Condition C2 in Theorem 1 contains that of the Relaxation conditions of group B in* [2].

*Proof*: The *Condition C2* given in *Theorem 1* is

$$\lambda_{j,t}^{\text{p}}>\lambda_{j,t}^{\text{p,C2}} \quad (68)$$

And the *condition C3_2* in *Relaxation conditions of group B* in [2] is

$$\lambda_{j,t}^{\text{p}}\geq \lambda_{j,t}^{\text{p,C3}} \quad (69)$$

Then, we can conclude that if $\lambda_{j,t}^{\text{p,C3}}>\lambda_{j,t}^{\text{p,C2}}$, *Lemma 4* holds. Below, we will provide the proof of $\lambda_{j,t}^{\text{p,C3}}>\lambda_{j,t}^{\text{p,C2}}$.

According to the *condition C3_1* and $0<\eta_{\text{ESS},n}^{\text{ch}}\cdot \eta_{\text{ESS},n}^{\text{dc}}<1$, we have

$$\frac{\partial \text{obj}}{\partial\left(p_{\text{ESS},n,t}^{\text{ch}}\right)}\cdot \eta_{\text{ESS},n}^{\text{ch}}\cdot \eta_{\text{ESS},n}^{\text{dc}}>-\frac{\partial \text{obj}}{\partial\left(p_{\text{ESS},n,t}^{\text{dc}}\right)}\cdot \eta_{\text{ESS},n}^{\text{ch}}\cdot \eta_{\text{ESS},n}^{\text{dc}} \quad (70)$$

Add term $-\frac{\partial \text{obj}}{\partial\left(p_{\text{ESS},n,t}^{\text{ch}}\right)}$ to both sides of formula (70), and it follows that

$$-\frac{\partial \text{obj}}{\partial\left(p_{\text{ESS},n,t}^{\text{ch}}\right)}+\frac{\partial \text{obj}}{\partial\left(p_{\text{ESS},n,t}^{\text{ch}}\right)}\cdot \eta_{\text{ESS},n}^{\text{ch}}\cdot \eta_{\text{ESS},n}^{\text{dc}}>$$
$$-\frac{\partial \text{obj}}{\partial\left(p_{\text{ESS},n,t}^{\text{ch}}\right)}-\frac{\partial \text{obj}}{\partial\left(p_{\text{ESS},n,t}^{\text{dc}}\right)}\cdot \eta_{\text{ESS},n}^{\text{ch}}\cdot \eta_{\text{ESS},n}^{\text{dc}} \quad (71)$$

Further, divide both sides of inequality (71) by $\left(1-\eta_{\text{ESS},n}^{\text{ch}} \cdot \eta_{\text{ESS},n}^{\text{dc}}\right)$, we can obtains

$$\lambda_{j,t}^{\text{p,C3}} = -\frac{\partial \text{obj}}{\partial \left(p_{\text{ESS},n,t}^{\text{ch}}\right)} >$$

$$\left[\frac{\partial \text{obj}}{\partial \left(p_{\text{ESS},n,t}^{\text{ch}}\right)} / \eta_{\text{ESS},n}^{\text{dc}} + \frac{\partial \text{obj}}{\partial \left(p_{\text{ESS},n,t}^{\text{dc}}\right)} \cdot \eta_{\text{ESS},n}^{\text{ch}}\right] / \left(\eta_{\text{ESS},n}^{\text{ch}} - \frac{1}{\eta_{\text{ESS},n}^{\text{dc}}}\right) \quad (72)$$

$$= \lambda_{j,t}^{\text{p,C2}}$$

Thus, **Lemma 4** holds. ∎

### B. Proof for Lemma 5

**Lemma 5**: *The region about LMPs formulated by Condition C2 in Theorem 1 contains that of the Relaxation conditions in* [5].

**Proof**: Recall that the condition given in **Condition C2** is

$$\lambda_{j,t}^{\text{p}} > \lambda_{j,t}^{\text{p,C2}} \quad (73)$$

Compared with **condition C7_1**, we can know that if $\lambda_{j,t}^{\text{p,C2}} < 0$, **Lemma 5** holds. Below, we will provide the proof of $\lambda_{j,t}^{\text{p,C2}} < 0$ on the premise that **condition C7** holds. We will give proofs by contradiction.

i) Assume $p_{\text{ESS},n,t}^{\text{ch}}, p_{\text{ESS},n,t}^{\text{dc}} > 0$, then we have $\lambda_{\text{ESS},n,t}^{\text{ch,1}} = \lambda_{\text{ESS},n,t}^{\text{dc,1}} = 0$. Further, with (37), it follows that

$$\Gamma_{\text{ESS},n,t} = -\frac{\partial \text{obj}}{\partial \left(p_{\text{ESS},n,t}^{\text{ch}}\right)} / \left(\eta_{\text{ESS},n}^{\text{ch}} \cdot \Delta t\right)$$

$$- \left(\lambda_{j,t}^{\text{p}} - \lambda_{\text{ESS},n,t}^{\text{ch,1}}\right) / \left(\eta_{\text{ESS},n}^{\text{ch}} \cdot \Delta t\right)$$

$$- \left(\lambda_{\text{ESS},n,t}^{\text{ch,2}} + \frac{\lambda_{\text{ESS},n,t}^{\text{relax}}}{\overline{P}_{\text{ESS},n}^{\text{ch}}} + 2\lambda_{\text{ESS},n,t}^{\text{S,2}} \cdot p_{\text{ESS},n,t}^{\text{ch}} + \lambda_{\text{ESS},n,t}^{\text{RD,1}}\right) / \left(\eta_{\text{ESS},n}^{\text{ch}} \cdot \Delta t\right) \quad (74)$$

Combining **condition C7** and $\lambda_{\text{ESS},n,t}^{\text{ch,1}} = 0$, $\eta_{\text{ESS},n}^{\text{ch}} > 0$, $\lambda_{\text{ESS},n,t}^{\text{ch,2}}, \lambda_{\text{ESS},n,t}^{\text{relax}}, \lambda_{\text{ESS},n,t}^{\text{S,2}}, \lambda_{\text{ESS},n,t}^{\text{RD,1}} \geq 0$, we can obtain formula (75) from (74).

$$\Gamma_{\text{ESS},n,t} \leq -\frac{\partial \text{obj}}{\partial \left(p_{\text{ESS},n,t}^{\text{ch}}\right)} / \left(\eta_{\text{ESS},n}^{\text{ch}} \cdot \Delta t\right) \quad (75)$$

ii) Subsequently, we can derive formula (76) by summing formulas (37) and (38).

$$0 = \lambda_{\text{ESS},n,t}^{\text{ch,1}} + \lambda_{\text{ESS},n,t}^{\text{dc,1}} = \frac{\partial \text{obj}}{\partial \left(p_{\text{ESS},n,t}^{\text{ch}}\right)} + \frac{\partial \text{obj}}{\partial \left(p_{\text{ESS},n,t}^{\text{dc}}\right)} + \frac{\lambda_{\text{ESS},n,t}^{\text{relax}}}{\overline{P}_{\text{ESS},n}^{\text{ch}}} + \frac{\lambda_{\text{ESS},n,t}^{\text{relax}}}{\overline{P}_{\text{ESS},n}^{\text{dc}}}$$

$$+ \lambda_{\text{ESS},n,t}^{\text{ch,2}} + \lambda_{\text{ESS},n,t}^{\text{dc,2}} + 2\left(\lambda_{\text{ESS},n,t}^{\text{S,2}} \cdot p_{\text{ESS},n,t}^{\text{ch}} + \lambda_{\text{ESS},n,t}^{\text{S,1}} \cdot p_{\text{ESS},n,t}^{\text{dc}}\right) \quad (76)$$

$$+ \lambda_{\text{ESS},n,t}^{\text{RD,1}} + \lambda_{\text{ESS},n,t}^{\text{RU,1}} - \left(1/\eta_{\text{ESS},n}^{\text{dc}} - \eta_{\text{ESS},n}^{\text{ch}}\right) \cdot \Gamma_{\text{ESS},n,t} \cdot \Delta t$$

Further, equation (76) yields:

$$\Gamma_{\text{ESS},n,t} = \left(\frac{\partial \text{obj}}{\partial \left(p_{\text{ESS},n,t}^{\text{ch}}\right)} + \frac{\partial \text{obj}}{\partial \left(p_{\text{ESS},n,t}^{\text{dc}}\right)}\right) / \left[\left(\frac{1}{\eta_{\text{ESS},n}^{\text{dc}}} - \eta_{\text{ESS},n}^{\text{ch}}\right) \cdot \Delta t\right]$$

$$+ \left(\frac{\lambda_{\text{ESS},n,t}^{\text{relax}}}{\overline{P}_{\text{ESS},n}^{\text{ch}}} + \frac{\lambda_{\text{ESS},n,t}^{\text{relax}}}{\overline{P}_{\text{ESS},n}^{\text{dc}}} + \lambda_{\text{ESS},n,t}^{\text{ch,2}} + \lambda_{\text{ESS},n,t}^{\text{dc,2}}\right) / \left[\left(\frac{1}{\eta_{\text{ESS},n}^{\text{dc}}} - \eta_{\text{ESS},n}^{\text{ch}}\right) \cdot \Delta t\right] \quad (77)$$

$$+ 2\left(\lambda_{\text{ESS},n,t}^{\text{S,2}} \cdot p_{\text{ESS},n,t}^{\text{ch}} + \lambda_{\text{ESS},n,t}^{\text{S,1}} \cdot p_{\text{ESS},n,t}^{\text{dc}}\right) / \left[\left(1/\eta_{\text{ESS},n}^{\text{dc}} - \eta_{\text{ESS},n}^{\text{ch}}\right) \cdot \Delta t\right]$$

$$+ \left(\lambda_{\text{ESS},n,t}^{\text{RD,1}} + \lambda_{\text{ESS},n,t}^{\text{RU,1}}\right) / \left[\left(1/\eta_{\text{ESS},n}^{\text{dc}} - \eta_{\text{ESS},n}^{\text{ch}}\right) \cdot \Delta t\right]$$

Given $\lambda_{\text{ESS},n,t}^{\text{ch,2}}, \lambda_{\text{ESS},n,t}^{\text{dc,2}}, \lambda_{\text{ESS},n,t}^{\text{relax}}, \lambda_{\text{ESS},n,t}^{\text{S,1}}, \lambda_{\text{ESS},n,t}^{\text{S,2}}, p_{\text{ESS},n,t}^{\text{ch}}, p_{\text{ESS},n,t}^{\text{dc}} \geq 0$, $\lambda_{\text{ESS},n,t}^{\text{RD,1}}, \lambda_{\text{ESS},n,t}^{\text{RU,1}} \geq 0$ and $1/\eta_{\text{ESS},n}^{\text{dc}} > \eta_{\text{ESS},n}^{\text{ch}}$, we can obtain inequality (78) based on (77).

$$\Gamma_{\text{ESS},n,t} \geq \left(\frac{\partial \text{obj}}{\partial \left(p_{\text{ESS},n,t}^{\text{ch}}\right)} + \frac{\partial \text{obj}}{\partial \left(p_{\text{ESS},n,t}^{\text{dc}}\right)}\right) / \left[\left(1/\eta_{\text{ESS},n}^{\text{dc}} - \eta_{\text{ESS},n}^{\text{ch}}\right) \cdot \Delta t\right] \quad (78)$$

iii) When formulas (75) and (78) cannot be established simultaneously, there will be a contradiction, indicating that the assumption $p_{\text{ESS},n,t}^{\text{ch}}, p_{\text{ESS},n,t}^{\text{dc}} > 0$ is not valid. In this case, the relaxation is exact for the given **condition C7**. Thus, we can derive formula (79), which is essentially the necessary condition for the **Relaxation conditions** in [5] to hold.

$$\left(\frac{\partial \text{obj}}{\partial \left(p_{\text{ESS},n,t}^{\text{ch}}\right)} + \frac{\partial \text{obj}}{\partial \left(p_{\text{ESS},n,t}^{\text{dc}}\right)}\right) / \left[\left(1/\eta_{\text{ESS},n}^{\text{dc}} - \eta_{\text{ESS},n}^{\text{ch}}\right) \cdot \Delta t\right]$$

$$> -\frac{\partial \text{obj}}{\partial \left(p_{\text{ESS},n,t}^{\text{ch}}\right)} / \left(\eta_{\text{ESS},n}^{\text{ch}} \cdot \Delta t\right) \quad (79)$$

Further, reformulating (79), we have

$$0 < \frac{\partial \text{obj}}{\partial \left(p_{\text{ESS},n,t}^{\text{ch}}\right)} / \eta_{\text{ESS},n}^{\text{dc}} + \frac{\partial \text{obj}}{\partial \left(p_{\text{ESS},n,t}^{\text{dc}}\right)} \cdot \eta_{\text{ESS},n}^{\text{ch}} \quad (80)$$

With (80) and $1/\eta_{\text{ESS},n}^{\text{dc}} > \eta_{\text{ESS},n}^{\text{ch}}$, we can obtain

$$0 > \left[\frac{\partial \text{obj}}{\partial \left(p_{\text{ESS},n,t}^{\text{ch}}\right)} / \eta_{\text{ESS},n}^{\text{dc}} + \frac{\partial \text{obj}}{\partial \left(p_{\text{ESS},n,t}^{\text{dc}}\right)} \cdot \eta_{\text{ESS},n}^{\text{ch}}\right] / \left(\eta_{\text{ESS},n}^{\text{ch}} - \frac{1}{\eta_{\text{ESS},n}^{\text{dc}}}\right) = \lambda_{j,t}^{\text{p,C2}} \quad (81)$$

Thus, **Lemma 5** holds. ∎

**Discussion**: Comparing the **condition C5** and **condition C7**, we can find that **condition C5_2** and **condition C7_1** are consistent. Moreover, we also notice that formula (80) is exactly the **condition C5_1**. Then, combined with the aforementioned conclusion that formula (80) is the necessary condition for the **Relaxation conditions** in [5] to hold, we can infer that **condition C5_1** can be obtained when the **Relaxation conditions** in [5] hold. Therefore, it can be concluded that the **condition C5**, i.e., *Relaxation conditions of Group C* in [2], can be derived based on the **Relaxation conditions** in [5].

### C. Proof for Lemma 6

**Lemma 6**: *The region about LMPs formulated by Condition C2 in Theorem 1 contains those of the Relaxation conditions in* [1] *and the Relaxation conditions of group A in* [2].

**Proof**: The sufficient condition provided in **Condition C2** is

$$\lambda_{j,t}^{\text{p}} > \lambda_{j,t}^{\text{p,C2}} \quad (82)$$

And the **condition C4_2** is

$$\lambda_{j,t}^{\mathrm{p}} > \lambda_{j,t}^{\mathrm{p,C4}} \tag{83}$$

Then, we can conclude that if $\lambda_{j,t}^{\mathrm{p,C4}} \geq \lambda_{j,t}^{\mathrm{p,C2}}$, **Lemma 6** holds. Below, we will provide the proof of $\lambda_{j,t}^{\mathrm{p,C4}} \geq \lambda_{j,t}^{\mathrm{p,C2}}$.

According to the **condition C4_1** and $0 < \eta_{\mathrm{ESS},n}^{\mathrm{ch}} \cdot \eta_{\mathrm{ESS},n}^{\mathrm{dc}} < 1$, we have

$$\frac{\partial \mathrm{obj}}{\partial \left( p_{\mathrm{ESS},n,t}^{\mathrm{ch}} \right)} \cdot \eta_{\mathrm{ESS},n}^{\mathrm{ch}} \cdot \eta_{\mathrm{ESS},n}^{\mathrm{dc}} \geq -\frac{\partial \mathrm{obj}}{\partial \left( p_{\mathrm{ESS},n,t}^{\mathrm{dc}} \right)} \cdot \eta_{\mathrm{ESS},n}^{\mathrm{ch}} \cdot \eta_{\mathrm{ESS},n}^{\mathrm{dc}} \tag{84}$$

Add term $-\dfrac{\partial \mathrm{obj}}{\partial \left( p_{\mathrm{ESS},n,t}^{\mathrm{ch}} \right)}$ to both sides of formula (84), it follows that

$$-\frac{\partial \mathrm{obj}}{\partial \left( p_{\mathrm{ESS},n,t}^{\mathrm{ch}} \right)} + \frac{\partial \mathrm{obj}}{\partial \left( p_{\mathrm{ESS},n,t}^{\mathrm{ch}} \right)} \cdot \eta_{\mathrm{ESS},n}^{\mathrm{ch}} \cdot \eta_{\mathrm{ESS},n}^{\mathrm{dc}} \geq \\ -\frac{\partial \mathrm{obj}}{\partial \left( p_{\mathrm{ESS},n,t}^{\mathrm{ch}} \right)} - \frac{\partial \mathrm{obj}}{\partial \left( p_{\mathrm{ESS},n,t}^{\mathrm{dc}} \right)} \cdot \eta_{\mathrm{ESS},n}^{\mathrm{ch}} \cdot \eta_{\mathrm{ESS},n}^{\mathrm{dc}} \tag{85}$$

Further, divide both sides of inequality (85) by $\left(1 - \eta_{\mathrm{ESS},n}^{\mathrm{ch}} \cdot \eta_{\mathrm{ESS},n}^{\mathrm{dc}}\right)$, we can obtains

$$\lambda_{j,t}^{\mathrm{p,C4}} = -\frac{\partial \mathrm{obj}}{\partial \left( p_{\mathrm{ESS},n,t}^{\mathrm{ch}} \right)} \geq \\ \left[ \frac{\partial \mathrm{obj}}{\partial \left( p_{\mathrm{ESS},n,t}^{\mathrm{ch}} \right)} / \eta_{\mathrm{ESS},n}^{\mathrm{dc}} + \frac{\partial \mathrm{obj}}{\partial \left( p_{\mathrm{ESS},n,t}^{\mathrm{dc}} \right)} \cdot \eta_{\mathrm{ESS},n}^{\mathrm{ch}} \right] / \left( \eta_{\mathrm{ESS},n}^{\mathrm{ch}} - \frac{1}{\eta_{\mathrm{ESS},n}^{\mathrm{dc}}} \right) \tag{86} \\ = \lambda_{j,t}^{\mathrm{p,C2}}$$

Thus, **Lemma 6** holds. ∎

### D. Proof for Lemma 7

**Lemma 7**: *The region about LMPs formulated by Condition C2 in Theorem 1 contains that of the Relaxation conditions of group C in* [2].

*Proof*: **Condition C5_2** in **Relaxation conditions of group C** in [2] is

$$\lambda_{j,t}^{\mathrm{p}} \geq 0 \tag{87}$$

And the sufficient condition provided in **Condition C2** is

$$\lambda_{j,t}^{\mathrm{p}} > \lambda_{j,t}^{\mathrm{p,C2}} \tag{88}$$

Obviously, if $0 > \lambda_{j,t}^{\mathrm{p,C2}}$, **Lemma 7** holds.

Combining **condition C5_1** and $1/\eta_{\mathrm{ESS},n}^{\mathrm{dc}} > \eta_{\mathrm{ESS},n}^{\mathrm{ch}}$, we have

$$0 > \left[ \frac{\partial \mathrm{obj}}{\partial \left( p_{\mathrm{ESS},n,t}^{\mathrm{ch}} \right)} / \eta_{\mathrm{ESS},n}^{\mathrm{dc}} + \frac{\partial \mathrm{obj}}{\partial \left( p_{\mathrm{ESS},n,t}^{\mathrm{dc}} \right)} \cdot \eta_{\mathrm{ESS},n}^{\mathrm{ch}} \right] / \left( \eta_{\mathrm{ESS},n}^{\mathrm{ch}} - \frac{1}{\eta_{\mathrm{ESS},n}^{\mathrm{dc}}} \right) = \lambda_{j,t}^{\mathrm{p,C2}} \tag{89}$$

Therefore, **Lemma 7** holds. ∎

### E. Proof for Lemma 8

**Lemma 8**: *The region about LMPs formulated by Condition C2 in Theorem 1 contains that of the Relaxation conditions in* [3].

*Proof*: Combined with the proof in subsection D in Appendix C, we can know that the regions formulated by **Relaxation conditions of group C** in [2] contains that of **Relaxation conditions** in [3]. Then, with the proof of **Lemma 7**, we can draw the conclusion that **Lemma 8** holds. ∎

### F. Proof for Lemma 9

**Lemma 9**: *The region about LMPs formulated by Condition C2 in Theorem 1 contains that of the Relaxation conditions in* [4].

*Proof*: Comparing **Condition C6**, that is, **Relaxation conditions** in [3] and **Condition C8**, we can easily conclude that the region formulated by **Relaxation conditions** in [3] contains that of **Relaxation conditions** in [4]. Thus, combining the proof of **Lemma 8**, we can conclude that **Lemma 9** holds. ∎

## APPENDIX E
## PROOFS FOR THE RELAXATION CONDITIONS OF ESSs UNDER MULTI-PERIOD ED SCENARIOS

First, the mathematical model for multi-period ED is presented as follows:

$$\min (1)$$
$$\text{s.t. } (15)\text{-}(22), (25)\text{-}(26), (28)\text{-}(32)$$

$$\sum_{j \in \mathcal{N}} \left( \sum_{g \in \Phi_{\mathrm{G},j}} p_{g,t}^{\mathrm{G}} + \sum_{n \in \Phi_{\mathrm{RG},j}} p_{n,t}^{\mathrm{RG}} + \sum_{n \in \Phi_{\mathrm{ESS},j}} \left( p_{\mathrm{ESS},n,t}^{\mathrm{dc}} - p_{\mathrm{ESS},n,t}^{\mathrm{ch}} \right) - p_{j,t}^{\mathrm{D}} \right) = 0, \; \left( \lambda_{j,t}^{\mathrm{p}} \right) \tag{90}$$

$$-\overline{P}_{\mathrm{L},j} \leq \sum_{i \in \mathcal{N}} SF_{j-i} \left( \sum_{g \in \Phi_{\mathrm{G},i}} p_{g,t}^{\mathrm{G}} + \sum_{n \in \Phi_{\mathrm{RG},i}} p_{n,t}^{\mathrm{RG}} + \sum_{n \in \Phi_{\mathrm{ESS},i}} \left( p_{\mathrm{ESS},n,t}^{\mathrm{dc}} - p_{\mathrm{ESS},n,t}^{\mathrm{ch}} \right) - p_{i,t}^{\mathrm{D}} \right) \leq \overline{P}_{\mathrm{L},j}$$
$$\left( \lambda_{\mathrm{L},j}^{1}(t), \; \lambda_{\mathrm{L},j}^{2}(t) \right) \tag{91}$$

where $SF_{j-i}$ denotes the shift factor for bus $i$ on line $j$.

Subsequently, let's proceed to prove **Condition C1$^{\mathrm{ED}}$** and **C2$^{\mathrm{ED}}$**.

**Proof of Condition C1$^{\mathrm{ED}}$**: Denoting the Lagrangian function of the relaxed ED model as $L_{\mathrm{RM}}^{\mathrm{ED}}$, and defining a new variable $\Gamma_{\mathrm{ESS},n,t} = \sum_{\tau=t}^{T} \left(1 - \delta_{\mathrm{ESS},n}\right)^{\tau-t} \left( \lambda_{\mathrm{ESS},n,\tau}^{\mathrm{SOC},2} - \lambda_{\mathrm{ESS},n,\tau}^{\mathrm{SOC},1} \right)$. Employing KKT conditions, the following equations are obtained:

$$\frac{\partial L_{\mathrm{RM}}}{\partial \left( p_{\mathrm{ESS},n,t}^{\mathrm{ch}} \right)} = \frac{\partial \mathrm{obj}}{\partial \left( p_{\mathrm{ESS},n,t}^{\mathrm{ch}} \right)} + \mathrm{LMP}_{j,t} - \lambda_{\mathrm{ESS},n,t}^{\mathrm{ch},1} + \lambda_{\mathrm{ESS},n,t}^{\mathrm{ch},2} + \frac{\lambda_{\mathrm{ESS},n,t}^{\mathrm{relax}}}{\overline{P}_{\mathrm{ESS},n}^{\mathrm{ch}}} \\ + \eta_{\mathrm{ESS},n}^{\mathrm{ch}} \cdot \Gamma_{\mathrm{ESS},n,t} \cdot \Delta t + \lambda_{\mathrm{ESS},n,t}^{\mathrm{RD},1} = 0 \tag{92}$$

$$\frac{\partial L_{\mathrm{RM}}}{\partial \left( p_{\mathrm{ESS},n,t}^{\mathrm{dc}} \right)} = \frac{\partial \mathrm{obj}}{\partial \left( p_{\mathrm{ESS},n,t}^{\mathrm{dc}} \right)} - \mathrm{LMP}_{j,t} - \lambda_{\mathrm{ESS},n,t}^{\mathrm{dc},1} + \lambda_{\mathrm{ESS},n,t}^{\mathrm{dc},2} + \frac{\lambda_{\mathrm{ESS},n,t}^{\mathrm{relax}}}{\overline{P}_{\mathrm{ESS},n}^{\mathrm{dc}}} \\ -\Gamma_{\mathrm{ESS},n,t} \cdot \Delta t / \eta_{\mathrm{ESS},n}^{\mathrm{dc}} + \lambda_{\mathrm{ESS},n,t}^{\mathrm{RU},1} = 0 \tag{93}$$

where $\mathrm{LMP}_{j,t} = \lambda_{j,t}^{\mathrm{p}} + \sum_{j} SF_{j-i} \left( \lambda_{\mathrm{L},j}^{1}(t) - \lambda_{\mathrm{L},j}^{2}(t) \right)$ represents the location marginal price (LMP).

With $(92) / \eta_{\mathrm{ESS},n}^{\mathrm{dc}} + (93) \cdot \eta_{\mathrm{ESS},n}^{\mathrm{ch}}$, we have

$$\left[\frac{\partial \text{obj}}{\partial\left(p_{\text{ESS},n,t}^{\text{ch}}\right)} + \text{LMP}_{j,t} - \lambda_{\text{ESS},n,t}^{\text{ch},1} + \lambda_{\text{ESS},n,t}^{\text{ch},2} + \frac{\lambda_{\text{ESS},n,t}^{\text{relax}}}{\overline{P}_{\text{ESS},n}^{\text{ch}}}\right.$$
$$\left. + \lambda_{\text{ESS},n,t}^{\text{RD},1}\right]/\eta_{\text{ESS},n}^{\text{dc}} + \left[\frac{\partial \text{obj}}{\partial\left(p_{\text{ESS},n,t}^{\text{dc}}\right)} - \text{LMP}_{j,t} - \lambda_{\text{ESS},n,t}^{\text{dc},1}\right. \tag{94}$$
$$\left. + \lambda_{\text{ESS},n,t}^{\text{dc},2} + \frac{\lambda_{\text{ESS},n,t}^{\text{relax}}}{\overline{P}_{\text{ESS},n}^{\text{dc}}} + \lambda_{\text{ESS},n,t}^{\text{RU},1}\right] \cdot \eta_{\text{ESS},n}^{\text{ch}} = 0$$

Further, equation (94) yields:

$$\text{LMP}_{j,t} = \left[\frac{\partial \text{obj}}{\partial\left(p_{\text{ESS},n,t}^{\text{ch}}\right)}/\eta_{\text{ESS},n}^{\text{dc}} + \frac{\partial \text{obj}}{\partial\left(p_{\text{ESS},n,t}^{\text{dc}}\right)} \cdot \eta_{\text{ESS},n}^{\text{ch}}\right]/\left(\eta_{\text{ESS},n}^{\text{ch}} - \frac{1}{\eta_{\text{ESS},n}^{\text{dc}}}\right)$$
$$- \left\{\left[\lambda_{\text{ESS},n,t}^{\text{ch},2} + \frac{\lambda_{\text{ESS},n,t}^{\text{relax}}}{\overline{P}_{\text{ESS},n}^{\text{ch}}} + \lambda_{\text{ESS},n,t}^{\text{RD},1}\right]/\eta_{\text{ESS},n}^{\text{dc}}\right. \tag{95}$$
$$\left. + \left[\lambda_{\text{ESS},n,t}^{\text{dc},2} + \frac{\lambda_{\text{ESS},n,t}^{\text{relax}}}{\overline{P}_{\text{ESS},n}^{\text{dc}}} + \lambda_{\text{ESS},n,t}^{\text{RU},1}\right] \cdot \eta_{\text{ESS},n}^{\text{ch}}\right\}/\left(\frac{1}{\eta_{\text{ESS},n}^{\text{dc}}} - \eta_{\text{ESS},n}^{\text{ch}}\right)$$
$$+ \left(\lambda_{\text{ESS},n,t}^{\text{ch},1}/\eta_{\text{ESS},n}^{\text{dc}} + \lambda_{\text{ESS},n,t}^{\text{dc},1} \cdot \eta_{\text{ESS},n}^{\text{ch}}\right)/\left(1/\eta_{\text{ESS},n}^{\text{dc}} - \eta_{\text{ESS},n}^{\text{ch}}\right)$$

If *Condition C1*$^{\text{ED}}$ holds, formula (96) can be obtained from (95).

$$\left(\lambda_{\text{ESS},n,t}^{\text{ch},1}/\eta_{\text{ESS},n}^{\text{dc}} + \lambda_{\text{ESS},n,t}^{\text{dc},1} \cdot \eta_{\text{ESS},n}^{\text{ch}}\right)/\left(1/\eta_{\text{ESS},n}^{\text{dc}} - \eta_{\text{ESS},n}^{\text{ch}}\right) > 0 \tag{96}$$

Since $0 < \eta_{\text{ESS},n}^{\text{dc}}, \eta_{\text{ESS},n}^{\text{ch}} < 1$, we have

$$\lambda_{\text{ESS},n,t}^{\text{ch},1}/\eta_{\text{ESS},n}^{\text{dc}} + \lambda_{\text{ESS},n,t}^{\text{dc},1} \cdot \eta_{\text{ESS},n}^{\text{ch}} > 0 \tag{97}$$

Based on (97), we can further derive inequality (98) in consideration of $\eta_{\text{ESS},n}^{\text{dc}}, \eta_{\text{ESS},n}^{\text{ch}} > 0$ and $\lambda_{\text{ESS},n,t}^{\text{ch},1}, \lambda_{\text{ESS},n,t}^{\text{dc},1} \geq 0$.

$$\lambda_{\text{ESS},n,t}^{\text{ch},1} + \lambda_{\text{ESS},n,t}^{\text{dc},1} > 0 \tag{98}$$

(98) indicates that at least one of $p_{\text{ESS},n,t}^{\text{ch}} = 0$ and $p_{\text{ESS},n,t}^{\text{dc}} = 0$ holds considering the complementary slackness conditions of constraints (15) and (16), which means that $p_{\text{ESS},n,t}^{\text{ch}} \cdot p_{\text{ESS},n,t}^{\text{dc}} = 0$. Hence, under this condition, even if the complementary constraints are relaxed, the optimal solution still implicitly avoids the SCD. ∎

**Proof of Condition C2**$^{\text{ED}}$: In combination with the expression of $\text{LMP}_{j,t}$ in equation (95), we can conclude that if *Condition C2*$^{\text{ED}}$ holds, we can obtain (99).

$$\left(\lambda_{\text{ESS},n,t}^{\text{ch},1}/\eta_{\text{ESS},n}^{\text{dc}} + \lambda_{\text{ESS},n,t}^{\text{dc},1} \cdot \eta_{\text{ESS},n}^{\text{ch}}\right)/\left(1/\eta_{\text{ESS},n}^{\text{dc}} - \eta_{\text{ESS},n}^{\text{ch}}\right) >$$
$$\left\{\left[\lambda_{\text{ESS},n,t}^{\text{ch},2} + \frac{\lambda_{\text{ESS},n,t}^{\text{relax}}}{\overline{P}_{\text{ESS},n}^{\text{ch}}} + \lambda_{\text{ESS},n,t}^{\text{RD},1}\right]/\eta_{\text{ESS},n}^{\text{dc}}\right. \tag{99}$$
$$\left. + \left[\lambda_{\text{ESS},n,t}^{\text{dc},2} + \frac{\lambda_{\text{ESS},n,t}^{\text{relax}}}{\overline{P}_{\text{ESS},n}^{\text{dc}}} + \lambda_{\text{ESS},n,t}^{\text{RU},1}\right] \cdot \eta_{\text{ESS},n}^{\text{ch}}\right\}/\left(\frac{1}{\eta_{\text{ESS},n}^{\text{dc}}} - \eta_{\text{ESS},n}^{\text{ch}}\right)$$

Considering $\eta_{\text{ESS},n}^{\text{dc}}, \eta_{\text{ESS},n}^{\text{ch}} > 0$, $\lambda_{\text{ESS},n,t}^{\text{ch},2}, \lambda_{\text{ESS},n,t}^{\text{dc},2} \geq 0$ and $\lambda_{\text{ESS},n,t}^{\text{relax}}, \lambda_{\text{ESS},n,t}^{\text{RD},1}, \lambda_{\text{ESS},n,t}^{\text{RU},1}, p_{\text{ESS},n,t}^{\text{ch}}, p_{\text{ESS},n,t}^{\text{dc}} \geq 0$, we can derive

$$\left[\lambda_{\text{ESS},n,t}^{\text{ch},2} + \frac{\lambda_{\text{ESS},n,t}^{\text{relax}}}{\overline{P}_{\text{ESS},n}^{\text{ch}}} + \lambda_{\text{ESS},n,t}^{\text{RD},1}\right]/\eta_{\text{ESS},n}^{\text{dc}}$$
$$+ \left[\lambda_{\text{ESS},n,t}^{\text{dc},2} + \frac{\lambda_{\text{ESS},n,t}^{\text{relax}}}{\overline{P}_{\text{ESS},n}^{\text{dc}}} + \lambda_{\text{ESS},n,t}^{\text{RU},1}\right] \cdot \eta_{\text{ESS},n}^{\text{ch}} \geq 0 \tag{100}$$

Combining formulas (99), (100) and $1/\eta_{\text{ESS},n}^{\text{dc}} > \eta_{\text{ESS},n}^{\text{ch}}$, we have

$$\left(\lambda_{\text{ESS},n,t}^{\text{ch},1}/\eta_{\text{ESS},n}^{\text{dc}} + \lambda_{\text{ESS},n,t}^{\text{dc},1} \cdot \eta_{\text{ESS},n}^{\text{ch}}\right) > 0 \tag{101}$$

Since $\eta_{\text{ESS},n}^{\text{dc}}, \eta_{\text{ESS},n}^{\text{ch}} > 0$ and $\lambda_{\text{ESS},n,t}^{\text{ch},1}, \lambda_{\text{ESS},n,t}^{\text{dc},1} \geq 0$, (102) holds.

$$\lambda_{\text{ESS},n,t}^{\text{ch},1} + \lambda_{\text{ESS},n,t}^{\text{dc},1} > 0 \tag{102}$$

Formula (102) indicates that $p_{\text{ESS},n,t}^{\text{ch}} \cdot p_{\text{ESS},n,t}^{\text{dc}} = 0$. In other words, based on *Condition C2*$^{\text{ED}}$, the optimal solution enforces the complementary constraint, even though it is relaxed. ∎